\documentclass[12pt]{article}
\usepackage{amsmath,latexsym,epsf}
\usepackage{amsfonts}
\usepackage{amssymb}
\usepackage{mathrsfs} 
\usepackage{amsmath}
\usepackage[english]{babel} 
\usepackage[latin1]{inputenc}
\usepackage[all]{xy}
\newcounter{Exe}[section]
\renewcommand{\theExe}{{\noindent \thesection.\arabic{Exe}}}
\newcommand{\Exe}{\vskip 10pt\noindent\refstepcounter{Exe}\textbf{%
           \theExe}\hskip 10pt}

\newcounter{Exesub}[Exe]
\renewcommand{\theExesub}{{\noindent \theExe.\arabic{Exesub}}}
\newcommand{\Exesub}{\vskip 10pt\noindent\refstepcounter{Exesub}\textbf{%
           \theExesub}\hskip 10pt}

\newcounter{Exesubsub}[Exesub]
\renewcommand{\theExesubsub}{{\noindent \theExesub\ (\roman{Exesubsub})}}
\newcommand{\Exesubsub}{\vskip 10pt\noindent\refstepcounter{Exesubsub}\textbf{%
           \theExesubsub}\hskip 10pt}

\begin{document}

\newcommand{\nc}{\newcommand}
\def\PP#1#2#3{{\mathrm{Pres}}^{#1}_{#2}{#3}\setcounter{equation}{0}}
\def\ns{$n$-star}\setcounter{equation}{0}
\def\nt{$n$-tilting}\setcounter{equation}{0}
\def\Ht#1#2#3{{{\mathrm{Hom}}_{#1}({#2},{#3})}\setcounter{equation}{0}}
\def\qp#1{{${(#1)}$-quasi-projective}\setcounter{equation}{0}}
\def\mr#1{{{\mathrm{#1}}}\setcounter{equation}{0}}
\def\mc#1{{{\mathcal{#1}}}\setcounter{equation}{0}}
\def\ms#1{{{\mathscr{#1}}}\setcounter{equation}{0}}
\def\mb#1{{{\mathbf{#1}}}\setcounter{equation}{0}}
\def\mf#1{{{\mathfrak{#1}}}\setcounter{equation}{0}}
\newcommand{\LL}{\ell\ell\,}
\newcommand{\edge}{\ar@{-}}
\newcommand{\bxy}{\xymatrix}
\newtheorem{Th}{Theorem}[section]
\newtheorem{Def}[Th]{Definition}
\newtheorem{Lem}[Th]{Lemma}
\newtheorem{Pro}[Th]{Proposition}
\newtheorem{Cor}[Th]{Corollary}
\newtheorem{Rem}[Th]{Remark}
\newtheorem{Exm}[Th]{Example}
\newtheorem{Obs}[Th]{Observation}
\def\Pf#1{{\noindent\bf Proof}.\setcounter{equation}{0}}
\def\>#1{{ $\Rightarrow$ }\setcounter{equation}{0}}
\def\<>#1{{ $\Leftrightarrow$ }\setcounter{equation}{0}}
\def\bskip#1{{ \vskip 20pt }\setcounter{equation}{0}}
\def\sskip#1{{ \vskip 5pt }\setcounter{equation}{0}}
\def\zskip#1{{ \vskip 10pt }\setcounter{equation}{0}}
\def\bg#1{\begin{#1}\setcounter{equation}{0}}
\def\ed#1{\end{#1}\setcounter{equation}{0}}
\def\r#1{{\rm{#1}}}
\def\adr{\mr{add}_{\mc{D}^b(\mr{mod}R)}}
\def\ads{\mr{add}_{\mc{D}^b(\mr{mod}S)}}
\def\kbr{\ms{K}^b(\mc{P}_R)}
\def\kbri{\ms{K}^b(\mc{I}_R)}
\def\kbs{\ms{K}^b(\mc{P}_S)}
\def\dbr{\ms{D}^b(R)}
\def\dbs{\ms{D}^b(S)}
\def\dr#1{\ms{D}^{#1}(R)}
\def\ds#1{\ms{D}^{#1}(S)}
\def\kfr{\ms{K}^{-,b}(\mc{P}_R)}
\def\kfri{\ms{K}^{+,b}(\mc{I}_R)}
\def\kfs{\ms{K}^{-,b}(\mc{P}_S)}
\def\ks#1{\ms{K}^{#1}(\mc{P}_S)}
\def\kr#1{\ms{K}^{#1}(\mc{P}_R)}
\def\kri#1{\ms{K}^{#1}(\mc{I}_R)}
\def\Ok#1{\Omega_{\ms{D}}^{#1}}
\def\Oki#1{\Omega^{\ms{D}}_{#1}}
\def\Hr{\mr{Hom}_{\ms{D}(R)}}
\def\Hs{\mr{Hom}_{\ms{D}(S)}}

\def\hjz#1#2{{ \begin{pmatrix} {#1} & {#2} \end{pmatrix}}\setcounter{equation}{0}}
\def\ljz#1#2{{  \begin{pmatrix} {#1} \\ {#2} \end{pmatrix}}\setcounter{equation}{0}}
\def\jz#1#2#3#4{{  \begin{pmatrix} {#1} & {#2} \\ {#3} & {#4} \end{pmatrix}}\setcounter{equation}{0}}
\def\23jz#1#2#3#4#5#6{{  \begin{pmatrix} {#1} & {#2} & {#3}\\ {#4} & {#5} & {#6} \end{pmatrix}}\setcounter{equation}{0}}
\def\32jz#1#2#3#4#5#6{{  \begin{pmatrix} {#1} & {#2} \\ {#3} & {#4} \\  {#5} & {#6} \end{pmatrix}}\setcounter{equation}{0}}




\title{\bf Repetitive equivalences and good Wakamatsu-tilting modules \thanks{Supported by the National Science Foundation of China (Grant No. 11371196) and the National Science Foundation for Distinguished Young Scholars of Jiangsu Province (Grant No. BK2012044) and
a project funded by the Priority Academic Program Development of Jiangsu Higher Education Institutions.}}
\smallskip
\author{{\small Jiaqun WEI}\\%
\small Institute of Mathematics, School of Mathematics Sciences \\
\small Nanjing Normal University, 
Nanjing 210023, P.R.China\\ \small Email:
weijiaqun@njnu.edu.cn
}
\date{}
\maketitle
\baselineskip 18pt
%
%
\begin{abstract}
\vskip 10pt%
Let $R$ be a ring and $T$ be a good Wakamatsu-tilting module with $S=\mr{End}_RT$. We prove that  $T$ induces an equivalence between stable categories of repetitive algebras $\hat{R}$ and $\hat{S}$.

%
\zskip\
%


\end{abstract}
%
\vskip 10pt
%

%
\section{\large Introduction}
%
%
\ \ \
Tilting theory plays an important role in the representation theory of artin algebras. The classical tilting modules were introduced in the early eighties by Brenner-Butler \cite{BB}, Bongartz \cite{Bgz} and Happel and Ringel \cite{HR}. Beginning with Miyashita \cite{My} and Happel \cite{HR}, the defining conditions for a classical tilting module were relaxed to tilting modules of arbitrary finite projective dimension, and further were relaxed to arbitrary rings and infinitely generated modules by many authors such as Colby and Fuller \cite{CF1}, Colpi
and Trlifaj \cite{CT}, Angeleri-H\"{u}gel and Coelho \cite{AhC}, Bazzoni \cite{Bz} etc..

One important result in tilting theory is the famous Brenner-Butler Theorem which shows that a tilting module induces some equivalences
between certain subcategories. In this sense, tilting theory may be viewed as a far-reaching way of generalization of
 Morita theory of equivalences between module categories. More interesting, when considering the derived category of an algebra, which contains module category of the algebra as a full subcategory,  Happel \cite{Hpb} and later   Cline, Parshall and Scott \cite{CPS} proved that a tilting module of finite projective dimension induces an equivalence between the bounded derived category of the ordinary algebra and the derived category of the endomorphism algebra of the tilting module. This leads to the study of Morita theory for derived categories, which were completely solved by Rickard \cite{Ric} through the notion of tilting complexes and by Keller \cite{Kel} through dg-categories.

A further generalization of tilting modules to tilting modules of possibly infinite projective
dimension was given by Wakamatsu \cite{Wk1}. Following \cite{GRS},  such tilting modules of possibly infinite projective
dimension are called Wakamatsu-tilting modules. It is known that Wakamatsu-tilting modules  also induce some equivalences between certain subcategories of module categories \cite{Wk2}. But Wakamatsu-tilting modules don't  induce derived equivalences in general.


However, we will show in this paper that Wakamatsu-tilting modules make more sense when we consider a more general category than the derived category of an algebra, namely, the stable module category of the repetitive algebra  of an algebras. To be compared, let us call the later category the repetitive category of the algebra. The repetitive category is a triangulated category. Moreover, by Happel's result \cite{Hpb}, for an artin algebra $R$, there is a fully faithful triangle embedding of the bounded derived category of $R$ into the repetitive category of $R$. Moreover, this embedding is an equivalence if and only if the global dimension of $R$ is finite.

We say that two algebras are repetitive equivalent if there is an equivalence between their repetitive categories.
It should be noted that repetitive equivalences are more general than derived equivalences. In fact, by results in \cite{As, Chqh, Ric} etc., if two algebras are derived  equivalent, then their repeptitive algebras are derived equivalent, and hence stably
equivalent. Thus derived equivalences always induce repetitive equivalences.


The following is our main theorem.

\vskip 10pt

\noindent{\bf Main Theorem}\ \  {\it Let $R$ be an artin algebra. If $T$ is a good Wakamatsu-tilting $R$-module with $S=\mr{End}(T_R)$, i.e., bimodules $_ST_R$ and $_RDT_S$ represent a cotorsion pair counter equivalence between a complete hereditary cotorsion piar $(\mc{B},\mc{A})$ in $\mr{mod}R$ and a complete hereditary cotorsion piar $(\mc{G},\mc{K})$ in $\mr{mod}S$, then $R$ and $S$ are repetitive equivalent. The equivalence restricts to the equivalence between $\mc{A}$ and $\mc{G}$.}

\vskip 10pt

We refer to subsection \ref{gwt} for more details on  good Wakamatsu-tilting modules. Note that that examples of good Wakamatsu-tilting modules contain tilting modules of finite projective dimension and cotilting modules of finite injective dimension.

Though a good Wakamatsu-tilting module induces a repetitive equivalence, unfortunately we can't say anything about whether or not the equivalence is a triangle equivalence now. However, if a repetitive equivalence is a triangle equivalence, we have the following result.

\vskip 10pt
\noindent{\bf Proposition} {\it Let $R$ and $S$ are artin algebras. Assume that there is a triangle equivalence between their repetitive categories and that this equivalence restricts to an equivalence between a covariantly finite coresolving subcategory $\mc{A}$ in $\mr{mod}R$ and a contravariantly finite resolving subcategory $\mc{G}$ in $\mr{mod}R$. Let $T$ be the preimage in $\mr{mod}R$ of $S$. Then $T$ is a good Wakamatsu-tilting $R$-module with $S\simeq\mr{End}(T_R)$.
}
\vskip 10pt

The paper is organized as follows. After the introduction, we provide basic knowledge on Wakamatsu-tilting modules and repetitive categories in Section 2. Then in Section 3 we introduce good Wakamatsu-tilting modules through cotorsion pair counter equivalences. Some properties and characterizations of good Wakamatsu-tilting modules are presented. Section 4  is devoted to the proof of the main theorem and the proposition in the introduction. Though the proof of the theorem is a little  complicated, the main idea is inspired by constructions in \cite[Lemma 4.1 in Chapter 3]{Hpb} and \cite[Section 1]{Wk2}. Finally, we provide some examples in the last section. In particular, it is shown that every Wakamatsu-tilting module over an algebra of finite representation type is a good Wakamatsu-tilting module and hence induces a repetitive equivalence.

%
%
\vskip 10pt
\noindent {\bf Conventions} Throughout this paper, we always work over artin algebras and right modules unless we claim otherewise. For an algebra $R$, we denote by $\mr{mod}R$ the category of all finitely generated $R$-modules, and by $\mr{proj}R$ (resp., $\mr{inj}R$) the category of finitele generated projective (resp., injective) $R$-modules. We denote the usual duality over an artin algebra $R$ by $D$.

For two functors $\mr{F}:\mc{A}\to\mc{B}$ and $\mr{G}:\mc{B}\to \mc{C}$, we use $\mr{GF}$ to denote their composition. While we use $f\circ g$, or simply just $fg$, to denote the composition of two homomorphisms $f:A\to B$ and $g:B\to C$.

Let $\mr{F}:\mc{A}\to\mc{B}$ be  a functor, we use $\mr{KerF}$ to denote the subcategory of $A\in\mc{A}$ such that $\mr{F}(A)=0$. Moreover, if $\mr{F}_i: \mr{A}\to \mc{B}$, $i\in I$, is a class of functors, we denote $\mr{KerF}_{I}=\bigcap_{i\in I}\mr{KerF}$. For instance, $\mr{KerExt}_R^{\ge 1}(T,-)$ is the subcategory of all $M\in\mr{mod}R$ such that $\mr{Ext}_R^i(T,M)=0$ for all $i\ge 1$.

We write the elements of direct sums as row vectors.


%

%
\bskip\
%

%
\vskip 30pt
\section{{\large  Wakamatsu-tilting modules and repetitive categories}}


%
%

\Exe\label{Wts}
{\bf Wakamatsu-tilting modules}

\zskip\

Recall that an $R$-module $T$ is {\it Wakamatsu-tilting} [\ref{Wk1}] provided that

\sskip\

(1) $\mr{End}_ST\simeq R$, where $S:=\mr{End}_RT$ and,

(2) $\mr{Ext}_R^i(T,T)=0=\mr{Ext}_S^i(T,T)=0$ for all $i>0$.

\sskip\

\noindent These two conditions are also equivalent to the following two conditions  [\ref{Wk1}].

\sskip\

(1) $\mr{Ext}_R^i(T,T)=0$ for all $i>0$ and,

\hangafter 1\hangindent 40pt
(2) There is an exact sequence $0\to R\to T_0\to T_1\to\cdots$, where $T_i\in \mr{add}_RT$ for all $i$, which stays exact after applying the functor
$\mr{Hom}_R(-,T)$.

\sskip\

Note that if $T$ is Wakamatsu-tilting and $S=\mr{End}(T_R)$, then $_S T$ is a Wakamatsu-tilting left $S$-module. In this case, we say that $T$ is a Wakamatsu-tilting
$S$-$R$-bimodule. It is easy to see that $DT$ is a Wakamatsu-tilting $R$-$S$-bimodule in the mean time.
\zskip\

\Exesub\label{AC}
{\it Auslander class and co-Auslander class}

\zskip\

Let $T$ be a Wakamatsu-tilting module with $S=\mr{End}(T_R)$. There are the following  two interesting classes associated with Wakamatsu-tilting modules.

The {\it Auslander class} in $\mr{mod}R$ with respect to the Wakamatsu-tilting  module $T_R$, denoted by $\mc{X}_{T}$, is defined as follows [\ref{AR}].

\vskip 5pt\hskip 40pt
\hangafter 1\hangindent 100pt
$\mc{X}_T:=\{M\in\mr{mod}R |$ there is an infinite exact sequence $0\to M\stackrel{f_0}{\longrightarrow} T_0\stackrel{f_1}{\longrightarrow} T_1\stackrel{f_2}{\longrightarrow} \cdots$ such that
$\mr{Im}f_i\in\mr{KerExt}_R^{\ge 1}(-,T)$ for each $i\ge 0 \}$.
\vskip 5pt

Obviously, it hold that $\mc{X}_T\subseteq \mr{KerExt}_R^{\ge 1}(-,T)$. Moreover, the two classes coincide with each other provided that $T$ is a cotilting $R$-module.

Dually, the {\it co-Auslander class} in $\mr{mod}R$ with respect to the Wakamatsu-tilting $R$-module $T$, denoted by ${_T\mc{X}}$, is defined as follows.

\vskip 5pt\hskip 40pt
\hangafter 1\hangindent 100pt
${_T\mc{X}}:=\{M\in\mr{mod}R |$ there is an infinite exact sequence $ \cdots \stackrel{f_2}{\longrightarrow} T_1\stackrel{f_1}{\longrightarrow} T_0\stackrel{f_0}{\longrightarrow} M \to 0$ such that
$\mr{Im}f_i\in\mr{KerExt}_R^{\ge 1}(T,-)$ for each $i\ge 0 \}$.
\vskip 5pt

Similarly, we have that $_T\mc{X}\subseteq \mr{KerExt}_R^{\ge 1}(T,-)$ and they coincide with each other provided that $T$ is a tilting $R$-module.

\zskip\ %
The following result collects some properties about Auslander class and co-Auslander class for a Wakamatsu-tilting module \cite{AR, MR, Wk2, Wk3}.

\zskip\

\noindent {\bf Proposition} 
{\it %
Let $T$ be a Waksmatsu-tilting $R$-module with $S=\mr{End}_RT$.%

\hangafter 1\hangindent 40pt
$(1)$ The Auslander class $\mc{X}_T$ is a resolving subcategory, i.e., it contains all projective $R$-modules and is closed under extensions, kernels of epimorphisms and direct summands.%

\hangafter 1\hangindent 40pt
$(2)$ The co-Auslander class $_T\mc{X}$ is  a coresolving subcategory, i.e., it contains all injective $R$-modules and is closed under extensions, cokernels of monomorphisms and direct summands.%

$(3)$ $\mr{KerExt}_R^1(\mc{X}_T,-)=\mr{KerExt}_R^{\ge 1}(\mc{X}_T,-)\subseteq {_T\mc{X}}$.%

$(4)$ $\mr{KerExt}_R^{1}(-,{_T\mc{X}})=\mr{KerExt}_R^{\ge 1}(-,{_T\mc{X}})\subseteq \mc{X}_T$.%

\hangafter 1\hangindent 39pt
$(5)$ The functors $\mr{Hom}_R(T,-)$ and $-\otimes_ST$ induce an equivalence between the co-Auslander class  $_T\mc{X}$ in $\mr{mod}R$ and the Auslander class
$\mc{X}_{DT}$ in $\mr{mod}S$. The  equivalence restricts to an  equivalence between the class  $\mr{KerExt}_R^{\ge 1}(\mc{X}_T,-)$ and the class
$\mr{KerExt}_S^{\ge 1}(-,{_{DT}\mc{X}})$.

}

\vskip 5pt
\Pf. (1) and (2) follows from \cite[Section 5]{AR}, see also \cite{MR}.

(3) and (4) follows from \cite[Lemma 1.4 and Proposition 1.6]{Wk3}.

(5) follows from  \cite[Proposition 2.14]{Wk2}.
\hfill $\Box$

\vskip 10pt
We remark that in case $T=R$, the class $\mc{X}_T=\mc{X}_R$ is just the class of all Gorenstein projective $R$-modules. Dually, in case $T=DR$, the class ${_T\mc{X}}={_{DR}\mc{X}}$ is just the class of all Gorenstein injective modules. We refer to \cite{EJb} for more on Gorenstein projective and Gorenstein injective modules.
\zskip\

\Exesub\label{cxt}
The following is a characterization of Auslander class and co-Auslander class, by \cite[Section 2]{Wk2} .

\vskip 10pt
\noindent {\bf Lemma} {\it Let $T$ be a Wakamatsu-tilting  $R$-module with $S=\mr{End}(T_R)$. Assume $X\in\mr{mod}R$ }

\hangafter 1\hangindent 20pt
\noindent{\it $(1)$ $X\in{_T\mc{X}}$ if and only if $X\in\mr{KerExt}_R^{>0}(T,-)$, $\mr{Hom}_R(T,X)\in\mr{KerTor}^S_{>0}(-,T)$ and $\mr{Hom}_R(T,X)\otimes_ST\simeq X$ canonically.}

\hangafter 1\hangindent 20pt
\noindent{\it $(2)$ $X\in{\mc{X}_T}$ if and only if $X\in\mr{KerExt}_R^{>0}(-,T)$, $\mr{Hom}_R(X, T)\in\mr{KerExt}_S^{>0}(-,T)$ and $\mr{Hom}_R(\mr{Hom}_R(X,T),T)\simeq A$ canonically.}

\sskip\
\zskip\

\Exesub\label{isom}
{\it Useful isomorphisms}

\vskip 10pt
Let $T$ be a Wakamatsu-tilting  $S$-$R$-bimodule. Then we have the following isomorphisms of bimodules.

\sskip\

 \centerline{$_SDS_S\simeq {_ST\otimes_R DT_S}$ and $_RDR_R\simeq {_RDT\otimes_S T_R}$.}

\sskip\

Given an adjoint pair $(\mr{F},\mr{G})$ of functors, we denote by $\mb{\Gamma}$ the natural adjoint isomorphism

\vskip 5pt
\centerline{$\mb{\Gamma}:\ \mr{Hom}(\mr{F}(-),-)\simeq \mr{Hom}(-,\mr{G}(-))$.}
\vskip 5pt

\noindent Moreover, for a homomorphism  $f: \mr{F}(X)\to Y$, we denote by $\mb{\Gamma}(f): X\to \mr{G}(Y)$ the image of $f$ under the isomorphism  $\mb{\Gamma}$.

In particular, associated with a  $S$-$R$-bimodule $T$, we have the following adjoint isomorphism

\vskip 5pt
\centerline{$\mb{\Gamma}^T:\ \mr{Hom}_R(-\otimes_ST,-)\simeq \mr{Hom}_S(-,\mr{Hom}_R(T,-))$.}
\vskip 5pt

\noindent We denote by $\eta^T$ and $\epsilon^T$ the unit and counit of this adjoint pair respectively, i.e.,

\vskip 5pt
\centerline{$\eta^T_X=\mb{\Gamma}^T(1_{X\otimes_ST}): X\to \mr{Hom}_R(T,X\otimes_ST)$ and}

\vskip 5pt
\centerline{$\epsilon^T_Y=(\mb{\Gamma}^T)^{-1}(1_{\mr{Hom}_R(T,Y)}): \mr{Hom}_R(T,Y)\otimes_ST\to Y$ }

\vskip 5pt
\noindent for $X\in\mr{mod}S$ and $Y\in\mr{mod}R$ respectively.

By the naturality of the isomorphism $\mb{\Gamma}$, for all homomorphisms $f:X_1\to X_2$, $g: \mr{F}(X_2)\to Y_1$ and $h:Y_1
\to Y_2$, it holds that $\mb{\Gamma}(\mr{F}(f)\circ g\circ h)=f\circ \mb{\Gamma}(g)\circ \mr{G}(h)$.

\zskip\

\Exe\label{rp}
{\bf Reprtitive algebras and Repetitive categories}
\vskip 10pt


\Exesub\label{rpa}
%
 We recall some basic facts on repetitive algebras mainly from \cite{Hpb}.

Let $R$ be an artin algebra. The repetitive algebra $\widehat{R}$ of $R$ was first introduced in \cite{HW}), which is defined to be the direct sum $\widehat{R}=\bigoplus_{n\in\mathbb{Z}}R\oplus\bigoplus_{n\in\mathbb{Z}}DR$ with the multiplication given by

\centerline{$(a_n,\varphi_n)(b_n,\psi_n)_n=(a_nb_n,a_{n+1}\psi_n+\varphi_nb_n)_n$.}

\noindent The repetitive algebra $\widehat{R}$ can be interpreted as the following infinite matrix algebra (without the identity)

\vskip 10pt
{\footnotesize\centerline{$\left(\begin{array}{ccccc}\ddots&&&&\\
\ddots & R & & & \\
& DR &R & &\\
&&DR&R&\\
&&&\ddots&\ddots
\end{array}\right)$.}
}

\vskip 10pt


\Exesub\label{rpm}
Denote by $\mr{mod}\widehat{R}$ the category of finitely generated $\widehat{R}$-modules. Then $\mr{mod}\widehat{R}$ is equivalent to the following two equivalent categories.

\vskip 5pt
\hangafter 1\hangindent 40pt
(1) $\mc{RC}^{\otimes}(R):=\{X=\{X_i,\delta^{\otimes}_i(X)\}_{i\in\mathbb{Z}}\mid X_i\in\mr{mod}R $ such that almost all $X_i$ are 0 and that $ \delta^{\otimes}_i(X):X_i\otimes_RDR\to X_{i-1}$ satisfying $(\delta^{\otimes}_{i+1}(X)\otimes_RDR)\circ\delta^{\otimes}_i(X)=0$, for each $i\}$.

\vskip 5pt
\hangafter 1\hangindent 40pt
(2) $\mc{RC}^{\mr{H}}(R):=\{X=\{X_i,\delta^{\mr{H}}_i(X)\}_{i\in\mathbb{Z}}\mid X_i\in\mr{mod}R $ such that almost all $X_i$ are 0 and that $\delta^{\mr{H}}_i(X):X_i\to \mr{Hom}_R(DR,X_{i-1})$ satisfying $\delta^{\mr{H}}_{i+1}(X)\circ  \mr{Hom}_R(DR,\delta^{\mr{H}}_i(X))=0$, for each $i\}$.

\vskip 5pt

We will freely use these equivalences. In particular, we often view $X\in\mr{mod}\widehat{R}$ as the following form with almost terms $X_i=0$

\vskip 5pt
\centerline{$ \cdots\stackrel{\delta_{i+1}}{\rightsquigarrow}X_i \stackrel{\delta_{i}}{\rightsquigarrow}X_{i-1}\stackrel{\delta_{i-1}}{\rightsquigarrow}\cdots$.}

\noindent and we call it a (bounded chain) repe-complex with the repe-difference $\delta$. We denote by $\mc{RC}(R)$ the category of all such repe-complexes. Thus, $\mc{RC}(R)=\mr{mod}\widehat{R}$. Note that there is an obvious automorphism $[1]:\mc{RC}(R)\to\mc{RC}(R)$ defined by $(X[1])_i=X_{i-1}$ for each $i$.

We say that a repe-complex $X=\{X_i,\delta_i\}\in\mc{RC}(R)$ is trivial  if each $\delta_i=0$. The full subcategory of all trivial repe-complexes is denoted by $\mc{RC}^{\mr{tr}}(R)$. Note that there is a natural forgetting functor from $\mc{RC}(R)$ to $\mc{RC}^{\mr{tr}}(R)$ by forgetting the repe-difference. 

Let $\mc{C}$ be a class of $R$-modules, we denote by $\mc{RC}(\mc{C})$ the class of repe-complexes with terms in $\mc{C}$. The notation $\mc{RC}^{\mr{tr}}(\mc{C})$ is defined similarly.

\zskip\

\Exesub\label{rpc}
%
As shown in \cite{Hpb}, $\widehat{R}$ is a selfinjective algebra and the category $\mc{RC}(R)(=\mr{mod}\widehat{R})$ is a Frobenius category, where the projective (and also injective) objects are of the form

\vskip 5pt

\centerline{$ \cdots\stackrel{\delta_{i+1}}{\rightsquigarrow}P_i\oplus I_i \stackrel{\delta_{i}}{\rightsquigarrow}P_{i-1}\oplus I_{i-1}\stackrel{\delta_{i-1}}{\rightsquigarrow}\cdots$,}

\vskip 5pt
\noindent where $P_i\in\mr{proj}R$, $I_i\in\mr{inj}R$ and $\delta_i=\jz{0}{\delta'_i}{0}{0}$ such that $\delta'_i: P_i\otimes_RDR\to I_{i-1}$ is an isomorphism (considered in $\mc{RC}^{\otimes}(R)$), or equivalently, $\delta'_i: P_i\to \mr{Hom}_R(DR, I_{i-1})$ is an isomorphism (considered in $\mc{RC}^{\mr{H}}(R)$).
Thus its stable category $\underline{\mc{RC}}(R)$ is a triangulated category. To compare with the derived category of an algebra, we will call it the {\it repetitive category} of $\mr{mod}R$ (or simply, of $R$).

It was shown in \cite{Hpb} that there is a fully faithful triangle embedding from the derived category $\mc{D}^b(\mr{mod}R)$ to the repetitive category $\underline{\mc{RC}}(R)$. Moreover, there is a triangle equivalence between $\mc{D}^b(\mr{mod}R)$ and $\underline{\mc{RC}}(R)$ if and only if $R$ has finite global dimension.

\vskip 5pt
For basic knowledge on triangulated categories,  derived
categories and the tilting theory, we refer to \cite{Hpb}.


\zskip\
\zskip\
\section{{\large Cotorsion pairs and good Wakamatusu-tilting modules}}

\Exe\label{cotp}
{\bf Cotorsion pair counter equivalences}

\vskip 10pt

A pair of subcategories $(\mc{B},\mc{A})$ in $\mr{mod}R$ is called a cotorsion pair, if $\mc{B}=\mr{KerExt}_R^1(-,\mc{A})$ and  $\mc{A}=\mr{KerExt}_R^1(\mc{B},-)$. A cotorsion pair $(\mc{B},\mc{A})$ is called hereditary provided that $\mc{B}$ is resolving, or equivalently, $\mc{A}$ is coresolving. Moreover, a cotorsion pair $(\mc{B},\mc{A})$ is called complete provided that, for each $X\in\mr{mod}R$, there exist exact sequences $0\to X\to A\to B\to 0$ and $0\to A'\to B'\to X\to 0$ for some $A,A'\in\mc{A}$ and $B,B'\in\mc{B}$. We refer to the book \cite{GT} for the general results on cotorsion pairs.

Let $(\mc{B},\mc{A})$ be a cotorsion pair in $\mr{mod}R$ and  $(\mc{G},\mc{K})$  be a cotorsion pair in $\mr{mod}S$. Similarly to torsion theory counter equivalences studied in \cite{CF, CDT},  we say that there is a {\it cotorsion pair counter equivalence} between $(\mc{B},\mc{A})$ and  $(\mc{G},\mc{K})$ provided that there are an equivalence $\mr{H}: \mc{A}\  {^{\longrightarrow}_{\longleftarrow}} \ \mc{G} :\mr{T}$ and  an equivalence $\mr{H}^{\prime}: \mc{K} \ {^{\longrightarrow}_{\longleftarrow}}\ \mc{B} :\mr{T}^{\prime}$. Moreover, we say that two bimodules $_SV_R$ and $_RV_S'$ represent the cotorsion pair counter equivalence if $\mr{H}=\mr{Hom}_R(V,-)$, $\mr{T}=-\otimes_SV$ and  $\mr{H}'=\mr{Hom}_S(V',-)$, $\mr{T}=-\otimes_RV'$.

There are close relations between Wakamatsu-tilting modules and cotorsion pairs, as shown in the following proposition.

\vskip 10pt

\noindent {\bf Proposition} {\it Let $T$ be a Waksmatsu-tilting $R$-module with $S=\mr{End}(T_R)$.} 

\hangafter 1\hangindent 39pt
{\it $(1)$ Both pairs $(\mr{KerExt}_R^1(-,{_T\mc{X}}),{_T\mc{X}})$ and $(\mc{X}_T,\mr{KerExt}_R^1(\mc{X}_T,-))$ are hereditary cotorsion pairs.}%

\hangafter 1\hangindent 39pt
{\it $(2)$ The bimodules $_ST_R$ and $_RDT_S$ represent a cotorsion pair counter equivalence between the cotorsion pair $(\mr{KerExt}_R^1(-,{_T\mc{X}}),{_T\mc{X}})$ in $\mr{mod}R$ and  the cotorsion pair $(\mc{X}_{DT},\mr{KerExt}_S^1(\mc{X}_{DT},-))$ in $\mr{mod}S$.}

\hangafter 1\hangindent 39pt
{\it $(3)$ The  bimodules $_RDT_S$ and $_ST_R$ represent a cotorsion pair counter equivalence between the cotorsion pair $(\mr{KerExt}_S^1(-,{_{DT}\mc{X}}),{_{DT}\mc{X}})$  in $\mr{mod}S$ and the cotorsion pair $(\mc{X}_{T},\mr{KerExt}_R^1(\mc{X}_{T},-))$  in $\mr{mod}R$.
}

%
%
\vskip 5pt
\hangafter 1\hangindent 40pt
\Pf. (1) follows from \cite[Proposition 3.1]{MR} and Propsition \ref{AC}. \\
(2) follows from Propositon \ref{AC} (5).\\
 (3) is obtained from (2) by replacing $_ST_R$ with $_RDT_S$.
%
\hfill $\Box$

\zskip\

\Exe\label{gwt}
{\bf Good Wakamatsu-tilting modules}

\vskip 10pt

\Exesub\label{dgwt}
%
%
In general case, the two cotorsion pairs in Proposition \ref{cotp} (1) are not complete. For instance, consider the case $T=R$. Then $\mc{X}_R$ is the class of all Gorensten projective modules. It is well known that this class is not a precovering class in general, see for instance \cite{Yo}. Thus, the cotorsion pair $(\mc{X}_R,\mr{KerExt}_R^1(\mc{X}_R,-))$ cannot be complete. Dually, in case $T=DR$, the cotorsion pair $(\mr{KerExt}_R^1(-,{_{DR}\mc{X}}),{_{DR}\mc{X}})$ is not complete in general.

However, the other cotorsion pair of  the  two cotorsion pairs in Proposition \ref{cotp} (1) for the above examples is complete respectively. This leads to the following general definition.
%
%
\vskip 10pt

\noindent {\bf Definition} {\it A Wakamatsu-tilting bimodule $_ST_R$ is said to be {\bf good} if the bimodules $_ST_R$ and $_RDT_S$ represent a cotorsion pair counter equivalence between a compete hereditary cotorsion pair $(\mc{B},\mc{A})$ in $\mr{mod}R$ and a compete hereditary cotorsion pair $(\mc{G},\mc{K})$ in $\mr{mod}S$. Furthermore, an $R$-module $T$ is said to be  a good Wakamatsu-tilting module if $_ST_R$ is a good Wakamatsu-tilting bimodule with $S=\mr{End}(T_R)$}.

\vskip 10pt
For example, $R$ and $DR$ are good Wakamatsu-tilting modules.

\vskip 10pt

\Exesub\label{pgwt}
By the definition, we have the following property about good Wakamatsu-tilting bimodules.

\vskip 10pt
\hangafter 1\hangindent 39pt
\noindent {\bf Proposition} {\it Let $_ST_R$ be a good Waksmatsu-tilting bimodule. Assume that  $(\mc{B},\mc{A})$ is  a compete hereditary cotorsion pair in $\mr{mod}R$  and  $(\mc{G},\mc{K})$ is a compete hereditary cotorsion pair in $\mr{mod}S$ such that the bimodules $_ST_R$ and $_RDT_S$ represent a cotorsion pair counter equivalence between them. Then

\hangafter 1\hangindent 39pt
$(1)$\ \ There is an equivalence \\
 \centerline{$\mr{Hom}_R(T,-):\hskip 20pt \mc{A}\ \ \  {^{\longrightarrow}_{\longleftarrow}}\ \ \ \mc{G}\hskip 20pt  :-\otimes_ST$\ \ \ }
  and  an equivalence \\
  \centerline{$\mr{Hom}_S(DT,-):\hskip 20pt \mc{K}\ \ \  {^{\longrightarrow}_{\longleftarrow}}\ \ \ \mc{B}\hskip 20pt  :-\otimes_RDT$.}

\hangafter 1\hangindent 39pt
$(2)$\ \  $_RDT_S$ is also a good Waksmatsu-tilting bimodule.

$(3)$\ \  $\mc{B}\subseteq\mc{X}_T$, $\mc{A}\in {_T\mc{X}}$ and $\mc{G}\subseteq\mc{X}_{DT}$, $\mc{K}\subseteq {_{DT}\mc{X}}$.

$(4)$\ \  $\mr{add}_RT= \mc{B}\bigcap\mc{A}$ and $\mr{add}_SDT= \mc{G}\bigcap\mc{K}$.
}

\vskip 10pt
\Pf. (1)  follows from the definition of good Wakamatsu-tilting bimodules.

(2) Replacing the Waksmatsu-tilting bimodule $_ST_R$ with the Wakamatsu-tilting bimodule $_RDT_S$ and noting that $DDT=T$, one can obtain (2) directly.

(3)\ \  Firstly, we show that   $\mr{add}_RT\subseteq\mc{B}\bigcap\mc{A}$ and $\mr{add}_SDT\subseteq \mc{G}\bigcap\mc{K}$.

 Note that all the involved subcategories in (1) are  closed under finite direct sums and direct summands.  Since $\mc{G}$ is resolving, we have that $S\in \mc{G}$. By the first equivalence in (1), we obtain that $T=S\otimes_ST\in\mc{A}$. It follows that $\mr{add}_RT\subseteq \mc{A}$. Dually, since $\mc{K}$ is coresolving, we have that $DS\in\mc{K}$. It follows from the second equivalence in (1) that $T=\mr{Hom}_S(S,T))=\mr{Hom}_S(DT,DS)\in\mc{B}$. Hence, $\mr{add}_RT\subseteq \mc{B}$ too. Thus, we obtain that  $\mr{add}_RT\subseteq\mc{B}\bigcap\mc{A}$.
Dually, one also has $\mr{add}_SDT\subseteq \mc{G}\bigcap\mc{K}$.

Clearly, $\mc{B}=\mr{KerExt}_S^1(-,\mc{A})=\mr{KerExt}_S^{\ge 1}(-,\mc{A})\subseteq\mr{KerExt}_S^{\ge 1}(-,T)$ follows from  $\mr{add}_RT\subseteq\mc{B}\bigcap\mc{A}$ and the fact that $(\mc{B},\mc{A})$ is a complete and hereditary cotorsion pair. Take any $B\in\mc{B}$, then $B\otimes_RDT\in\mc{K}$. Take an exact sequence $0\to B\otimes_RDT\to I\to Y\to 0$ with $I\in\mr{inj}S=\mr{add}_SDS$. Since $\mc{K}$ is coresolving, we have that $I,Y\in\mc{K}$ too. Applying the functor $\mr{Hom}_S(DT,-)$, we obtain an induced exact sequence $0\to \mr{Hom}_S(DT,B\otimes_RDT)\to \mr{Hom}_S(DT,I)\to \mr{Hom}_S(DT,Y)\to 0$, since $\mc{K}=\mr{KerExt}_S^1(\mc{G},-)\subseteq\mr{KerExt}_S^1(DT,-)$ by the above argument. Note that $B\simeq \mr{Hom}_S(DT,B\otimes_RDT)$, $\mr{Hom}_S(DT,I)\in\mr{add}_RT$ and $\mr{Hom}_S(DT,Y)\in\mc{B}$, so one can easily see that $B\in\mc{X}_T$. Thus $\mc{B}\subseteq\mc{X}_T$. By the equivalence in Proposition \ref{cotp} (3), we also obtain that $\mc{K}\subseteq {_{DT}\mc{X}}$.

Now consider the good Wakamatsu-tilting module $_RDT_S$ and apply the above result, we can obtain that $\mc{G}\in\mc{X}_{DT}$ and that $\mc{A}\in {_T\mc{X}}$.

 %
(4)\ \ If $X\in\mc{B}\bigcap\mc{A}$, then $X\in\mc{B}$. Following from the proof of (3), we obtain that there is an exact sequence $0\to X\to T_X\to X'\to 0$ with $T_X\in\mr{add}_RT$ and $X'\in\mc{B}$. Since $X\in \mc{A}$ too, we have that $\mr{Ext}_R^1(X',X)=0$. It follows that the exact sequence splits. Hence $X\in\mr{add}_RT$. Together with the first claim in the proof of (3), we obtain that $\mr{add}_RT= \mc{B}\bigcap\mc{A}$. Dually, we also have that $\mr{add}_SDT= \mc{G}\bigcap\mc{K}$.
 \hfill $\Box$

%
%

%

\zskip\

\Exesub\label{homlem}
%
Recall that a subcategory $\mc{A}\in\mr{mod}R$ is covariantly finite (or, a preenveloping calss) if for any $X\in\mr{mod}R$, there is an object $A_X\in\mc{A}$ and a homomorphism $u_{_X}: X\to A_X$ such that $\mr{Hom}_R(u_{_X},A)$ is surjective for any object $A\in\mc{A}$, see for instance \cite{AR}. Dually,  a subcategory $\mc{B}\in\mr{mod}R$ is contravariantly finite (or, a precovering calss)  if for any $X\in\mr{mod}R$, there is an object $B_X\in\mc{A}$ and a homomorphism $v_X: B_X\to X$ such that $\mr{Hom}_R(B,v_X)$ is surjective for any object $B\in\mc{B}$.

Let $\mc{A}$ be a subcategory of $\mr{mod}R$. An $R$-module $T$ is said to be Ext-projective if $T\in\mc{A}\bigcap\mr{KerExt}_R^1(-,\mc{A})$. Moreover, it is said to be an Ext-projective generator if, for any $A\in\mc{A}$, there exists an exact sequence $0\to A'\to T_A\to A\to 0$ with $T_A\in\mr{add}_RT$ and $A'\in\mc{A}$.  Dually, an $R$-module $T$ is said to be an Ext-injective cogenerator if $T\in\mc{A}\bigcap\mr{KerExt}_R^1(\mc{A},-)$ and that,  for any $A\in\mc{A}$, there exists an exact sequence $0\to A\to T_A\to A'\to 0$ with $T_A\in\mr{add}_RT$ and $A'\in\mc{A}$.

\vskip 10pt
\noindent{\bf Lemma} {\it Let $\mc{A}$ be a subcategory closed under extensions and direct summands.}

\hangafter 1\hangindent 35pt
{\it $(1)$ Assume that $\mc{A}$ has an Ext-projective generator $T$. If $0\to X\to Y\to Z\to 0$ is an exact sequence which stays exact after the functor $\mr{Hom}_R(T,-)$, where $Y,Z\in\mc{A}$, then $X\in\mc{A}$ too.}

\hangafter 1\hangindent 35pt
{\it $(2)$  Assume that $\mc{A}$ has an Ext-injective cogenerator $T$. If $0\to X\to Y\to Z\to 0$ is an exact sequence which stays exact after the functor $\mr{Hom}_R(-,T)$, where $X,Y\in\mc{A}$, then $Z\in\mc{A}$ too.
}

\vskip 5pt
\noindent\Pf. (1) By the assumptions, we can construct the following commutative diagram, where $T_Z\in\mr{add}_RT$ and $Z'\in\mc{A}$.

 \hskip 120pt$\xymatrix@C=20pt@R=20pt{
               &             & 0\ar[d]        &       0\ar[d]          & \\
               &             & Z'\ar[d]\ar@{=}[r]&      Z'\ar[d]   & \\
          0\ar[r]& X\ar[r]^{(1,0)\ \ \ }\ar@{=}[d] &X\oplus T_Z\ar[r]^{\ \ \ (^0_1)}\ar[d]^{(^f_h)} &T_Z\ar[dl]^h\ar[r]\ar[d]^t&0\\
          0\ar[r]& X\ar[r]_{f}  &Y\ar[r]_{g}\ar[d]  &Z\ar[r]\ar[d]&0\\
               &             & 0         &       0         &} $

\vskip 5pt
Since $\mc{A}$ is closed under extensions and direct summands, we have that $X\in\mc{A}$ from the middle column.

(2) Dually.
\hfill $\Box$

\vskip 10pt
%

\Exesub\label{res-cor}
\hangafter 1\hangindent 35pt
\noindent{\bf Lemma} {\it Let $T$ be a Wakamatsu-tilting $R$-module with $S=\mr{End}(T_R)$. Assume that $\mr{Hom}_R(T,-): \mc{A}\  {^{\longrightarrow}_{\longleftarrow}} \ \mc{G} : -\otimes_ST$ define an equivalence. Then the following are equivalent.}

{\it $(1)$ $\mc{A}$ is coresolving and $T$ is an Ext-projective generator in $\mc{A}$.}

{\it $(2)$ $\mc{G}$ is resolving and $DT$ is an Ext-injective cogenerator in $\mc{G}$.}

\vskip 5pt
\noindent \Pf. (1) $\Rightarrow$ (2) The condition that  $T$ is an Ext-projective generator in $\mc{A}$ says that $T\in\mc{A}\subseteq\mr{KerExt}_R^1(T,-)$ and that every $A\in\mc{A}$ admits an exact sequence $0\to A'\to T_A\to A\to 0$ with $T_A\in\mr{add}_RT$ and $A'\in\mc{A}$. This implies that $\mc{A}\subseteq\mr{KerExt}_R^{\ge 1}(T,-)$,  since $\mc{A}$ is coresolving. In particular, $\mc{A}\subseteq {_T\mc{X}}$.

Note that, for any $X\in\mc{A}$, there is an exact sequence $0\to X\to I\to X'\to 0$ with $I'\in\mr{inj}R\subseteq\mc{A}$ and $X'\in\mc{A}$. Applying the functor $\mr{Hom}_R(T,-)$, we have an exact sequence $0\to \mr{Hom}_R(T,X)\to \mr{Hom}_R(T,I)\to \mr{Hom}_R(T,X')\to 0$. Since $\mr{Hom}_R(T,I)\in\mr{add}_SDT$ and $\mr{Ext}_S^1(\mr{Hom}_R(T,X),DT)\simeq \mr{Ext}_S^1(\mr{Hom}_R(T,X),\mr{Hom}_R(T,DR))=0$,  we obtain that $DT$ is an Ext-injective cogenerator in $\mr{Hom}_R(T,\mc{A})=\mc{G}$.

It is clear that $\mc{G}$ is closed under direct summands. Assume now there is an exact sequence $(\flat):\  0\to X\to Y\stackrel{g}{\longrightarrow} Z\to 0$ with $Z\in\mc{G}$, then $Z\in\mr{Hom}_R(T,\mc{A})\subseteq\mr{KerTor}^S_1(-,T)$. Applying the functor $-\otimes_ST$, we obtain an induce exact sequence $(\flat\otimes_ST):\  0\to X\otimes_ST\to Y\otimes_ST\stackrel{g\otimes_ST}{\longrightarrow} Z\otimes_ST\to 0$. 

Assume first $X\in\mc{G}$ too, then $X\otimes_ST\in \mc{A}$. It follows that $Y\otimes_ST\in\mc{A}$ too, since $\mc{A}$ is closed under extensions. Note now there is an exact sequence  $0\to \mr{Hom}_R(T,X\otimes_ST)\to \mr{Hom}_R(T,Y\otimes_ST)\to \mr{Hom}_R(T,Z\otimes_ST)\to 0$, so we have that $\mr{Hom}_R(T,Y\otimes_ST)\simeq Y$ since $\mr{Hom}_R(T,X\otimes_ST)\simeq X$ and $\mr{Hom}_R(T,Z\otimes_ST)\simeq Z$. Thus $Y\in\mr{Hom}_R(T,\mc{A})=\mc{G}$. This shows that $\mc{G}$ is closed under extensions.

Assume now $Y\in\mc{G}$, then $\mr{Hom}_R(T,g\otimes_ST)\simeq g$. In particular, we have that  $\mr{Hom}_R(T,X\otimes_ST)\simeq X$ and the homomorphism $\mr{Hom}_R(T,g\otimes_ST)$ is surjective. It follows that the exact sequence $(\flat\otimes_ST)$ stays exact after the functor $\mr{Hom}_R(T,-)$. By Lemma \ref{homlem}, we obtain that $X\otimes_ST\in\mc{A}$. Hence $X\in \mr{Hom}_R(T,\mc{A})=\mc{G}$. This shows that $\mc{G}$ is closed under  kernels of epimorphisms.

(2) Dually.
%
\hfill $\Box$

\zskip\

\Exesub\label{cequ}
\noindent{\bf Proposition} {\it Let $_ST_R$ be a Wakamatsu-tilting module. Assume that $(\mc{B},\mc{A})$ is a hereditary cotorsion pair in $\mr{mod}R$ and that $T$ is an Ext-projective generator in $\mc{A}$, then $(\mr{Hom}_R(T,\mc{A}),\mc{B}\otimes_RDT)$ is a hereditary cotorsion pair in $\mr{mod}S$. In particular, the bimodules $_ST_R$ and $_RDT_S$ represent a cotorsion pair counter equivalence between $(\mc{B},\mc{A})$ and  $(\mr{Hom}_R(T,\mc{A}),\mc{B}\otimes_RDT)$ in the case.
}

\vskip 5pt
\noindent\Pf. Since $T$ is an Ext-projective generator in $\mc{A}$, we see that $T\in\mr{KerExt}_R^1(-,\mc{A})\bigcap \mc{A}=\mc{B}\bigcap \mc{A}$ and that, for any $A\in\mc{A}$, there is an exact sequence $0\to A'\to T_A\to A\to 0$ with $T_A\in\mr{add}_RT$ and $A'\in\mc{A}$. In particular, for any $X\in\mc{A}\bigcap\mc{B}$, there is an exact sequence $0\to X'\to T_X\to X\to 0$ with $T_X\in\mr{add}_RT$ and $X'\in\mc{A}$, which is clearly split. Hence, $X\in\mr{add}_RT$. It follows that $\mr{add}_RT=\mc{B}\bigcap \mc{A}$. Moreover, by  an argument similar to the one used in the proof of \cite[Proposition 2.13(b)]{MR}, we have that $T$ is also an Ext-injective cogenerator in $\mc{B}$. Note that these  imply that $\mc{A}\subseteq {_T\mc{X}}$ and that  $\mc{B}\subseteq {\mc{X}_T}$.

By Lemma \ref{res-cor}, we see that  $\mr{Hom}_R(T,\mc{A})$ is resolving and that $\mc{B}\otimes_RDT$ is coresolving. It is also clear that the bimodules $_ST_R$ and $_RDT_S$ represent a counter equivalence between two pairs $(\mc{B},\mc{A})$ and  $(\mr{Hom}_R(T,\mc{A}),\mc{B}\otimes_RDT)$, by assumptions.  So, it just remains to show that $(\mr{Hom}_R(T,\mc{A}),\mc{B}\otimes_RDT)$ is a cotorsion pair.

We divide the remained proof into three steps.

\vskip 5pt
\hangafter 1\hangindent 40pt
\noindent {\bf Step 1}  $\mr{Ext}^i_S(\mr{Hom}_R(T,A),B\otimes_RDT)=0$, for any $A\in \mc{A}$ and any $B\in\mc{B}$ and for any $i\ge 0$.

\vskip 5pt
 Note that there is an isomorphism $D\mr{Hom}_S(S,B\otimes_SDT)\simeq \mr{Hom}_R(B,S\otimes_ST)$ and that it induces an isomorphism $D\mr{Hom}_S(S_i,B\otimes_SDT)\simeq \mr{Hom}_R(B,S_i\otimes_ST)$, for any $S_i\in\mr{add}_SS$.

Now take $A\in\mc{A}$, since $T$ is an Ext-projective generator, there is a long exact sequence

\hskip 90pt{$\cdots\to T_n\to\cdots \to T_1\to T_0\to A\to 0$}\hfill $(\dag)$

\noindent with each $T_i\in\mr{add}_RT$ and each image in $\mc{A}$. Since  $B\in\mr{KerExt}_R^1(-,A)$, we have an induced exact sequence $\mr{Hom}_R(B,\dag)$:

\hskip 20pt{\footnotesize$\cdots\to \mr{Hom}_R(B,T_n)\to\cdots \to \mr{Hom}_R(B,T_1)\to  \mr{Hom}_R(B,T_0)\to  \mr{Hom}_R(B,A)\to 0$.}

 On the other hand, by applying the functor $D\mr{Hom}_S(\mr{Hom}_R(T,-),B\otimes_RDT)$, we have a complex $D\mr{Hom}_S(\mr{Hom}_R(T,\dag),B\otimes_RDT)$:

\vskip 10pt

{\footnotesize\hskip 20pt{$\cdots\to D\mr{Hom}_S(\mr{Hom}_R(T,T_n),B\otimes_RDT)\to\cdots \to D\mr{Hom}_S(\mr{Hom}_R(T,T_1),B\otimes_RDT)$}

 \hskip 40pt{$\ \ \ \ \to  D\mr{Hom}_S(\mr{Hom}_R(T,T_0),B\otimes_RDT)\to  D\mr{Hom}_S(\mr{Hom}_R(T,A),B\otimes_RDT)\to 0$.}}

Since $\mr{Hom}_R(T,\dag)$ is exact, one sees that the functor  $D\mr{Hom}_S(\mr{Hom}_R(T,-),B\otimes_RDT)$ is right exact. By the above isomorphism, we obtain the following isomorphisms of complexes

\hskip 40pt{ {\footnotesize$D\mr{Hom}_S(\mr{Hom}_R(T,\dag),B\otimes_RDT)\simeq \mr{Hom}_R(B,\mr{Hom}_R(T,\dag)\otimes_ST)\simeq \mr{Hom}_R(B,\dag)$.}}

\noindent But the later is exact, so we obtain that

\vskip 10pt
\hskip 40pt{\footnotesize{$\mr{Ext}^i_S(\mr{Hom}_R(T,A),B\otimes_RDT)\simeq \mr{H}^i(\mr{Hom}_S(\mr{Hom}_R(T,\dag),B\otimes_RDT))$

\hskip 40pt $\simeq D\mr{H}^i(D\mr{Hom}_S(\mr{Hom}_R(T,\dag),B\otimes_RDT))\simeq D\mr{Hom}_R(B,\dag)=0$.}}

Thus, step 1 is established. In particular, we obtain that $\mr{Hom}_R(T,\mc{A})
\subseteq\mr{KerExt}^1_S(-,\mc{B}\otimes_RDT)$ and that $\mr{Hom}_R(T,\mc{A})
\subseteq\mr{KerExt}^1_S(-,\mc{B}\otimes_RDT)$, due to the arbitrarity of $A\in\mc{A}$ and $B\in\mc{B}$.

\vskip 5pt
\hangafter 1\hangindent 40pt
\noindent {\bf Step 2}\ \  $\mr{KerExt}^1_S(-,\mc{B}\otimes_RDT)\subseteq \mr{Hom}_R(T,\mc{A})$.

\vskip 5pt

Take any $Y\in \mr{KerExt}^1_S(-,\mc{B}\otimes_RDT)$ and a projective resolution of $Y$:

\vskip 5pt
\hskip 90pt $\cdots\stackrel{f_{n+1}}{\longrightarrow} S_n\stackrel{f_n}{\longrightarrow}\cdots\stackrel{f_1}{\longrightarrow} S_1\stackrel{f_1}{\longrightarrow} S_0\stackrel{f_0}{\longrightarrow} Y\to 0$. \hfill ($\sharp$)

Note that $DT=R\otimes_RDT\in\mc{B}\otimes_RDT$ and that  $\mc{B}\otimes_RDT $ is coresolving, so we obtain that

\hskip 10pt{\footnotesize $Y\in\mr{KerExt}^1_S(-,\mc{B}\otimes_RDT)=\mr{KerExt}_S^{>0}(-,\mc{B}\otimes_RDT)\subseteq\mr{KerExt}_S^{>0}(-,DT)=\mr{KerTor}^S_{>0}(-,T)$.}

\noindent Then we have an induced exact sequence:

\hskip 40pt $\cdots\to S_n\otimes_ST\to\cdots\to S_1\otimes_ST\to S_0\otimes_ST\to Y\otimes_ST\to 0$.  \hfill $(\sharp\otimes_ST)$

For any $B\in\mc{B}$, applying the left exact functor $\mr{Hom}_R(B,-)$, we obtain a complex $ \mr{Hom}_R(B,\sharp\otimes_ST)$:

\vskip 10pt
\hskip 40pt {\footnotesize $\cdots\to \mr{Hom}_R(B,S_n\otimes_ST)\to\cdots\to \mr{Hom}_R(B,S_1\otimes_ST)$

\hskip 90pt $\to \mr{Hom}_R(B,S_0\otimes_ST)\to \mr{Hom}_R(B,Y\otimes_ST)\to 0$.}

Applying the right exact functor $D\mr{Hom}_S(-,B\otimes_SDT)$ to the sequence ($\sharp$), we obtain a complex $D\mr{Hom}_S(\sharp,B\otimes_SDT)$:

\vskip 10pt
\hskip 40pt {\footnotesize $\cdots\to D\mr{Hom}_S(S_n,B\otimes_SDT)\to\cdots\to D\mr{Hom}_S(S_1,B\otimes_SDT)$

\hskip 90pt $\to D\mr{Hom}_S(S_0,B\otimes_SDT)\to D\mr{Hom}_S(Y,B\otimes_SDT)\to 0$.}

By the isomorphism in Step 1 again, we have isomorphisms of complexes

\hskip 90pt $D\mr{Hom}_S(\sharp,B\otimes_SDT)\simeq \mr{Hom}_R(B,\sharp\otimes_ST)$.

\noindent It follows that, for any $i\ge 2$,

\vskip 10pt
\hskip 70pt {\footnotesize $\mr{Ext}_R^1(B,Y_{i}\otimes_ST)\simeq\mr{H}^{i-1}(\mr{Hom}_R(B,\sharp\otimes_ST))$

\hskip 100pt $\simeq\mr{H}^{i-1}(D\mr{Hom}_S(\sharp,B\otimes_SDT))=D\mr{Ext}^{i-1}_R(Y,B\otimes_ST)=0$,}

\noindent where $Y_i=\mr{Im}f_i$, since $B\in\mr{KerExt}_R^1(-,T)$ and $\mr{Hom}_R(B,-)$ is left right exact. This shows that $Y_i\otimes_ST\in\mr{KerExt}_R^1(\mc{B},-)=\mc{A}$, for all $i\ge 2$. As $S\otimes_ST\simeq T\in\mc{A}$ and $\mc{A}$ is coresolving, we obtain that $Y\otimes_ST\in \mc{A}$ too. Then we have an induced exact sequence $\mr{Hom}_R(T,\sharp\otimes_ST)$:

\vskip 10pt
\hskip 40pt {\footnotesize $\cdots\to \mr{Hom}_R(T,S_n\otimes_ST)\to\cdots\to \mr{Hom}_R(T,S_1\otimes_ST)$

\hskip 90pt $\to \mr{Hom}_R(T,S_0\otimes_ST)\to \mr{Hom}_R(T,Y\otimes_ST)\to 0$.}

\noindent It follows that $Y\simeq \mr{Hom}_R(T,Y\otimes_ST)\in\mr{Hom}_R(T,\mc{A})$, since $S_i\simeq \mr{Hom}_R(T,S_i\otimes_ST)$ for each $i$. This shows that $\mr{KerExt}^1_S(-,\mc{B}\otimes_RDT)\subseteq\mr{Hom}_R(T,\mc{A})$. Together with Step 1, we obtain that $\mr{KerExt}^1_S(-,\mc{B}\otimes_RDT)=\mr{Hom}_R(T,\mc{A})$.

\vskip 5pt
\hangafter 1\hangindent 40pt
\noindent {\bf Step 3}\ \  $\mr{KerExt}^1_S(\mr{Hom}_R(T,\mc{A}),-)\subseteq \mc{B}\otimes_RDT$.

\vskip 5pt

Note that there is an isomorphism \\
\centerline{$\mr{Hom}_S(\mr{Hom}_R(T,A),DS)\simeq D\mr{Hom}_R(\mr{Hom}_S(DT,DS),A)$,}\\
for any $A\in\mr{mod}R$ and that it induces an isomorphism

\hskip 90pt $\mr{Hom}_S(\mr{Hom}_R(T,A),I_i)\simeq D\mr{Hom}_R(\mr{Hom}_S(DT,I_i),A)$, 

\noindent for any $I_i\in\mr{add}DS$.

Now take any $X\in\mr{KerExt}_S^1(\mr{Hom}_R(T,\mc{A}),-)$ and 
consider an injective resolution of $X$:

\vskip 5pt
\hskip 90pt $0\to X\stackrel{g_0}{\longrightarrow} I_0\stackrel{g_1}{\longrightarrow} I_1\stackrel{g_1}{\longrightarrow}\cdots\stackrel{g_{n-1}}{\longrightarrow} I_n \stackrel{g_n}{\longrightarrow} \cdots$. \hfill $(\natural)$

\noindent Since $DT\simeq \mr{Hom}_R(T,DR)\in\mr{Hom}_R(T,\mc{A})$ and $\mr{Hom}_R(T,\mc{A})$ is resolving, we obtain that

\hskip 20pt {\footnotesize $X\in \mr{KerExt}_S^1(\mr{Hom}_R(T,\mc{A}),-)=\mr{KerExt}_S^{>0}(\mr{Hom}_R(T,\mc{A}),-)\subseteq \mr{KerExt}_S^{>0}(DT,-)$}.

Thus, for any $A\in\mc{A}$, applying the functor $\mr{Hom}_S(\mr{Hom}_R(T,A),-)$, we have an induce exact complex $\mr{Hom}_S(\mr{Hom}_R(T,A),\natural)$:

\vskip 5pt
\hskip 80pt {\footnotesize $0\to \mr{Hom}_S(\mr{Hom}_R(T,A),X){\longrightarrow} \mr{Hom}_S(\mr{Hom}_R(T,A),I_0) $}

\hskip 110pt {\footnotesize ${\longrightarrow}\mr{Hom}_S(\mr{Hom}_R(T,A),I_1){\longrightarrow}\cdots {\longrightarrow} \mr{Hom}_S(\mr{Hom}_R(T,A),I_n) {\longrightarrow} \cdots$}.

On the other hand, we also have the following induced complex  $D\mr{Hom}_R(\mr{Hom}_S(DT,\natural),A)$, by applying the functor $D\mr{Hom}_R(\mr{Hom}_S(DT,-),A)$: 

\vskip 5pt
\hskip 80pt {\footnotesize $0\to D\mr{Hom}_R(\mr{Hom}_S(DT,X),A){\longrightarrow} D\mr{Hom}_R(\mr{Hom}_S(DT,I_0),A) $}

\hskip 110pt {\footnotesize ${\longrightarrow}D\mr{Hom}_R(\mr{Hom}_S(DT,I_1),A){\longrightarrow}\cdots {\longrightarrow} D\mr{Hom}_R(\mr{Hom}_S(DT,I_n),A) {\longrightarrow} \cdots$}.


By the above mentioned isomorphism $(\ast)$ and the fact that $\mr{Hom}_S(DT,\natural)$ is exact, we have that, for any $i\ge 1$,

\vskip 5pt
\hskip 70pt {\footnotesize $\mr{Ext}_R^1(\mr{Hom}_S(DT,X_{i+2}),A)\simeq \mr{H}^{i}(\mr{Hom}_R(\mr{Hom}_S(DT,\natural),A))$}

\hskip 100pt {\footnotesize   $\simeq D\mr{H}^{i}(D\mr{Hom}_R(\mr{Hom}_S(DT,\natural),A))\simeq D\mr{H}^{i}(\mr{Hom}_S(\mr{Hom}_R(T,A),\natural))=0$,}

\noindent where $X_i=\mr{Im}g_i$, since $A\in\mr{KerExt}_R^1(-,T)$ and $\mr{Hom}_R(-,A)$ is left exact. It follows that $\mr{Hom}_S(DT,X_{i+2})\in\mr{KerExt}_R^1(-,\mc{A})=\mc{B}$ for any $i\ge 1$. Since $\mc{B}$ is resolving, we also obtain each $\mr{Hom}_S(DT,X_{i})\in \mc{B}$ for for each $0\le i\le 2$, where $X_0:=X$. Then we have an induced exact sequence $\mr{Hom}_S(DT,\natural)\otimes_RDT$,  since $\mc{B}\subseteq\mr{KerExt}_R^1(-,T)=\mr{KerTor}^R_1(-,DT)$:

\hskip 40pt {\footnotesize $0\to \mr{Hom}_S(DT,X)\otimes_SDT\to \mr{Hom}_S(DT,I_0)\otimes_SDT$}

\hskip 80pt {\footnotesize $\to\mr{Hom}_S(DT,I_1)\otimes_SDT\to\cdots\to \mr{Hom}_S(DT,I_n)\otimes_SDT\to\cdots$.}

\noindent It follows that $X\simeq \mr{Hom}_S(DT,X)\otimes_SDT\in\mc{B}\otimes_SDT$, since $I_i\simeq \mr{Hom}_S(DT,I_i)\otimes_SDT$ for each $i$. This shows  that $\mr{KerExt}^1_S(\mr{Hom}_R(T,\mc{A}),-)\subseteq\mc{B}\otimes_RDT$. Together with Step 1, we obtain that $\mr{KerExt}^1_S(\mr{Hom}_R(T,\mc{A}),-)=\mc{B}\otimes_RDT$.

\vskip 5pt
Altogether, we obtain that $(\mr{Hom}_R(T,\mc{A}),\mc{B}\otimes_RDT)$ is a hereditary cotorsion pair.
\hfill $\Box$

\vskip 10pt


\Exesub\label{c-gwt}
\noindent{\bf Corollary} {\it Let $T$ be a Wakamatsu-tilting $R$-module with $S=\mr{End}(T_R)$. Assume that the functor $\mr{Hom}_R(T,-)$ gives an equivalence between a covariantly finite coresolving subcategory $\mc{A}$ in $\mr{mod}R$ and  a contravariantly finite resolving subcategory $\mc{G}$ in $\mr{mod}S$. If $T$ is an Ext-projective generator in $\mc{A}$, then $T$ is a good Wakamatsu-tilting module.}

\vskip 5pt
\noindent\Pf. Since $\mc{G}$ is  a contravariantly finite resolving subcategory in $\mr{mod}S$, there is a cotorsion pair $(\mc{G},\mr{KerExt}_R^1(\mc{G},-))$ in $\mr{mod}S$, by \cite[Proposition 1.10]{AR}. Dually,  there is a cotorsion pair $(\mr{KerExt}_R^1(-,\mc{A}),\mc{A})$ in $\mr{mod}R$, since $\mc{A}$ is  a covariantly finite coresolving subcategory $\mc{A}$ in $\mr{mod}R$. Note that both  cotrosion pairs are complete and hereditary. Since  $T$ is an Ext-projective generator in $\mc{A}$, by Proposition \ref{cequ}, the bimodules $_ST_R$ and $_RDT_S$ represent a cotorsion pair counter equivalence between the above two cotorsion pairs. Hence $T$ is a good Wakamatsu-tilting module by the definition.
%
%
%
\hfill $\Box$

%
\vskip 30pt
\section{\large The proof of the main results}

\ \ \ \ \

The whole section will be devoted to the proof of the two results mentioned in the introduction.

Let $R$ be an artin algebra and $T$ be a good Wakamatsu-tilting module with $S=\mr{End}(T_R)$. Then $_ST_R$ be a good Wakamatsu-tilting bimodule.  Assume that  $(\mc{B},\mc{A})$ is  a compete hereditary cotorsion pair in $\mr{mod}R$  and  $(\mc{G},\mc{K})$ is a compete hereditary cotorsion pair in $\mr{mod}S$ such that the bimodules $_ST_R$ and $_RDT_S$ represent a cotorsion pair counter equivalence between these two  cotorsion pairs.

The sketch of our proof on the main theorem is as follows.

Firstly we construct a functor  $\mr{L}_T: \mc{RC}^{\mr{tr}}(R)\to\mc{RC}(S)$ and a functor $-\hat{\otimes} DT :\mc{RC}(R)\to\mc{RC}(S)$. Then we give a natural homomorphism

\vskip 5pt
\centerline{$l^X_Y:\mr{Hom}_{\mc{RC}^{\mr{tr}}(R)}(X,Y)\to \mr{Hom}_{\mc{RC}(S)}(X\hat{\otimes}DT,\mr{L}_T(Y))$}
\vskip 5pt

\noindent which is functional in both variables. After this, associated with an object $X\in\mc{RC}(R)$, we use the condition that  $(\mc{B},\mc{A})$ is  a compete cotorsion pair in $\mr{mod}R$ to obtain an object $A_X\in\mc{RC}^{\mr{tr}}(\mc{A})$ and  establish a homomorphism $u_{_X}\in \mr{Hom}_{\mc{RC}^{\mr{tr}}(R)}(X,A_X)$. We then show that the assignment that setting $X\longmapsto \mr{Cok}(l(u_{_X}))$ induces our desired functor $\mb{S}_T: \underline{\mc{RC}}(R)\to\underline{\mc{RC}}(S)$. We use the dual method to construct another desired functor $\mb{Q}_{DT}:  \underline{\mc{RC}}(S)\to \underline{\mc{RC}}(R)$. Then we prove that there are natural  isomorphisms  $\mb{Q}_{DT}\mb{S}_T\simeq 1_{\underline{\mc{RC}}(R)}$ and $\mb{S}_T\mb{Q}_{DT}\simeq 1_{\underline{\mc{RC}}(R)}$.

\vskip 10pt

\Exe\label{rpc}\label{S}
{\bf From $\mc{RC}(R)$ to $\mc{RC}(S)$: the functor $\mb{S}_T$}

\Exesub\label{0-l}\label{FcL}
{\it The functor $\mr{L}_T: \mc{RC}^{\mr{tr}}(R)\to\mc{RC}(S)$
}

\vskip 10pt

Let $X=\{X_i\}\in\mc{RC}^{\mr{tr}}(R)$. We define $\mr{L}_T(X)\in\mc{RC}(S)$ as follows.
\sskip\
%
%
%

$(l1)$ The underlying module $\mr{L}_T(X)_i=\mr{Hom}_R(T,X_{i-1})\oplus X_i\otimes_R DT$ and,

\hangafter 1 \hangindent 41pt
$(l2)$ The structure map  $\delta^{\otimes}_i(\mr{L}_T(X)): \mr{L}_T(X)_i\to \mr{L}_T(X)_{i-1}$ is given by $\begin{pmatrix} 0 & \delta_{L_i}\\ 0 & 0\end{pmatrix}$,  where $\delta_{L_i}$ is the composition:
%
%

{\footnotesize \hskip 5pt\centerline{$\mr{Hom}_R(T,X_{i-1})\otimes_S DS \stackrel{\simeq}{\longrightarrow} \mr{Hom}_R(T,X_{i-1})\otimes_S T\otimes_R \mr{D}T\stackrel{\epsilon_{X_{i-1}}^T\otimes_R
\mr{D}T}{\longrightarrow} X_{i-1}\otimes_R DT.$}}
\vskip 5pt

From the functor  property of $\mr{Hom}_R(T,-)$ and $-\otimes_R DT$, one can easily see that $\mr{L}_T$ is a functor from $\mc{RC}^{\mr{tr}}(R)$ to $\mc{RC}(S)$. 

\vskip 10pt

\noindent{\bf Remark}  (1) If $X\in\mc{RC}^{\mr{tr}}(\mr{add}_RT)$, i.e.,  $X=\{X_i\}$ with each  $X_i\in\mr{add}_RT$, then $\mr{Hom}_R(T,X_{i-1})\in \mr{add}_SS$ and $\delta_{L_i}$ defined above is an isomorphism for each
$i$. It follows that $\mr{L}_T(X)$ is a projective object in $\mc{RC}(S)$ in the case.

(2) In particular, in case $T=R$, we obtain the functor $\mr{L}_R: \mc{RC}^{\mr{tr}}(R)\to\mc{RC}(R)$ which specially send objects in $\mc{RC}^{\mr{tr}}(\mr{add}_RR)$ to  a projective object in $\mc{RC}(R)$.

\zskip\

\Exesub\label{r-r}
{\it The functor $-\hat{\otimes} DT :\mc{RC}(R)\to\mc{RC}(S)$ 
 }
\vskip 10pt

Let $Y=\{Y_i,\delta^{\otimes}_i(Y)\}\in\mc{RC}(R)$. We  define $Y\hat{\otimes} DT\in\mc{RC}(S)$ by setting 

\sskip\

$(t1)$ the underlying module is $(Y\hat{\otimes} DT)_i=Y_i\otimes_RDT$ and,

$(t2)$ the structure map $\delta^{\otimes}_i(Y\hat{\otimes} DT)$ is given by the composition

\vskip 5pt
\hskip 25pt {\footnotesize $Y_i\otimes_RDT\otimes_SDS \stackrel{\simeq}{\longrightarrow} Y_i\otimes_RDT\otimes_ST\otimes_RDT  \stackrel{\simeq}{\longrightarrow}
Y_i\otimes_RDR\otimes_RDT\stackrel{\delta^{\otimes}_i(Y)\otimes_RDT}{\longrightarrow} Y_{i-1}\otimes_RDT$}.
\sskip\

From the functor property of $-\otimes_R DT$, one can see that $-\hat{\otimes} DT$ is a functor from $\mc{RC}(R)$ to $\mc{RC}(S)$.

\zskip\

\Exesub\label{lpro}\label{lhom}
{\it The homomorphism  $l^X_Y:\mr{Hom}_{\mc{RC}^{\mr{tr}}(R)}(X,Y)\to \mr{Hom}_{\mc{RC}(S)}(X\hat{\otimes}DT,\mr{L}_T(Y))$ 
}
%
\vskip 10pt
Recall that we have a forgotting functor from $\mc{RC}(R)$ to $\mc{RC}^{\mr{tr}}(R)$.
 For any $X\in\mc{RC}(R)$ and $Y\in\mc{RC}^{\mr{tr}}(R)$, there is a canonical homomorphism

\vskip 5pt
 \centerline{$l^X_Y: \mr{Hom}_{\mc{RC}^{\mr{tr}}(R)}(X,Y)\to \mr{Hom}_{\mc{RC}(S)}(X\hat{\otimes}DT,\mr{L}_T(Y))$ }
\vskip 5pt

 \noindent which is functional in both variables, defined by


\vskip 5pt
\centerline{$l^X_Y: u=\{u_i\}$ \hskip 10pt $\longmapsto$ \hskip 10pt $f=\{f_i\}$ with $f_i=\hjz{-\theta_{l_i},}{{u_i\otimes_RDT}}$,}
\vskip 5pt


\noindent where $\theta_{l_i}$ is given by the composition

\sskip\

\hskip 10pt {\footnotesize $X_i\otimes_RDT\stackrel{\eta^T_{X_i}}{\longrightarrow} \mr{Hom}_R(T,X_i\otimes_R DT\otimes_ST)\stackrel{\simeq}{\longrightarrow}\mr{Hom}_R(T,X_i\otimes_RDR)$

\sskip\
 \hskip 150pt $\stackrel{\mr{Hom}_R(T,\delta_{X_i})}{\longrightarrow} \mr{Hom}_R(T,
X_{i-1})\stackrel{\mr{Hom}_R(T,u_{i-1})}{\longrightarrow}\mr{Hom}_R(T,Y_{i-1})$.}

\vskip 10pt
\noindent{\bf Remark} Using the fact ${_RDR_R}\simeq {_R(DT\otimes_ST)_R}$ and the adjoint isomorphism \\
\centerline{$\mb{\Gamma}^{T}: \mr{Hom}_R(X_i\otimes_RDR,Y_{i-1})\simeq\mr{Hom}_S(X_i\otimes_RDT,\mr{Hom}_R(T,Y_{i-1}))$,}\\ one can easily check that $\theta_{l_i}$ is just the image of the natural homomorphism $\delta^{\otimes}_i(X) \circ u_{i-1}$ under $\Gamma^T$, i.e., $\theta_{l_i}=\mb{\Gamma}^{T}(\delta^{\otimes}_i(X) \circ u_{i-1})$.


%
\zskip\

It is easy to see that, for any commutative diagram in $\mc{RC}(R)$

\vskip 5pt
  \centerline{  $\xymatrix{X\ar[r]^{u}\ar[d]^{x} & Y\ar@<-1ex>[d]^{y}\ \\
                X'\ar[r]^{u'}         & Y',
        }$ }

\vskip 5pt
\noindent there is an induced commutative diagram in $\mc{RC}(S)$

\vskip 5pt
 \centerline{ $\xymatrix{X\hat{\otimes}DT\ar[r]^{l(u)}\ar[d]^{x\hat{\otimes}DT} & \mr{L}_T(Y)\ar@<-1ex>[d]^{y}\ \\
              X'\hat{\otimes}DT\ar[r]^{l(u')}                          & \mr{L}_T(Y').
        }$}

\vskip 5pt

\bskip\

%
%
%
\zskip\


\Exesub\label{resolution}
{\it A monomorphism $u_{_X}: X\to A_X$ with $A_X\in\mc{RC}^{\mr{tr}}(\mc{A})$, for $X\in\mc{RC}(R)$}.

%
%
\zskip\

%
%
%

Let $X=\{X_i\}\in\mc{RC}(R)$. Since $(\mc{B},\mc{A})$ is a complete cotorsion pair in $\mr{mod}R$, there are exact sequences $0\to X_i\stackrel{u_{_{X_i}}}{\longrightarrow} A_{X_i}\to B_{X_i}\to 0$ with $A_{X_i}\in \mc{A}$ and $B_{X_i}\in \mc{B}$, for each $i$. This gives an exact sequence $0\to X\stackrel{u_{_X}}{\longrightarrow} A_X\stackrel{\pi_{u_{_X}}}{\longrightarrow} B_X\to 0$ in $\mc{RC}^{\mr{tr}}(R)$ with $A_X=\{A_{X_i}\}\in\mc{RC}^{\mr{tr}}(\mc{A})$ and $B_X=\{B_{X_i}\}\in\mc{RC}^{\mr{tr}}(\mc{B})$.

Now let $Y=\{Y_i\}\in \mc{RC}(R)$ and  $h=\{h_i\}\in\mr{Hom}_{\mc{RC}(R)}(X,Y)$. Then we have an exact sequence $0\to Y\stackrel{u_{_Y}}{\longrightarrow} A_Y\stackrel{\pi_{u_{_Y}}}{\longrightarrow} B_Y\to 0$  in $\mc{RC}^{\mr{tr}}(R)$ with $A_Y=\{A_{Y_i}\}\in\mc{RC}^{\mr{tr}}(\mc{A})$ and $B_Y=\{B_{Y_i}\}\in\mc{RC}^{\mr{tr}}(\mc{B})$, as above.
Note that $\mc{B}\subseteq\mr{KerExt}_R^1(-,\mc{A})$,
it is easy to see that there is a homomorphism $h_{A}\in \mr{Hom}_{\mc{RC}^{\mr{tr}}(R)}(A_{X},A_{Y})$ and further $h_{B}\in \mr{Hom}_{\mc{RC}^{\mr{tr}}(R)}(B_{X},B_{Y})$,  such that the following diagram in $\mc{RC}^{\mr{tr}}(R)$ is commutative with exact rows.

\vskip 5pt
  \centerline{  $\xymatrix{
 0\ar[r] & X\ar[r]^{u_{_X}}\ar[d]^{h} & A_X\ar[r]^{\pi_{{_X}}}\ar[d]^{h_A} & B_X\ar[r]\ar[d]^{h_B}&0 \\
 0\ar[r] & Y\ar[r]^{u_{_Y}}           & A_Y\ar[r]^{\pi_{{_Y}}}             & B_X\ar[r] & 0 \\
        }$ }

\bskip\

\Exesub\label{coker}
{\it The cokernel $\mr{Cok}(l(u_{_X}))$
}
\zskip\

Applying the functor $-\otimes_RDT$ to the exact sequences $0\to X_i\stackrel{u_{_{X_i}}}{\longrightarrow} A_{X_i}\to B_{X_i}\to 0$ in the above, we obtain an induced exact sequences

\vskip 5pt 
  \centerline{  $\xymatrix{
 0\ar[r] & X_i\otimes_RDT \ar@{->}[rr]^{u_{_{X_i}}\otimes_RDT} & & A_{X_i}\otimes_RDT\ar[r] & B_{X_i}\otimes_RDT\ar[r] & 0,}$
 }

\vskip 5pt

 \noindent
  since $B_{X_i}\in\mc{B}\subseteq\mr{KerExt}_R^1(-,T)=\mr{KerTor}^R_1(-,DT)$ for each $i$. It follows that, by applying the homomorphism $l$ in \ref{lpro} to the homomorphism $u_{_X}$ in \ref{resolution}, there is an induced exact sequence

\centerline{$\xymatrix{
0\ar[r] & X\hat{\otimes}DT\ar@{->}[rr]^{l(u_{_X})}  & & \mr{L}_T(A_X)\ar@{->}[rr]^{\pi_{_{l_X}}} & & \mr{Cok}(l_{_X})\ar[r]  & 0.} $ }

\vskip 10pt

\noindent {\bf Remark} From the definition of $\mr{Cok}(l(u_{_X}))$, one see that, for each $i$, $\mr{Cok}(l(u_{_X}))_i$ is given by the pushout

  \centerline{  $\xymatrix{
             {X_i}\otimes_RDT\ar[r]^{\theta_{l_i}\  \  }\ar[d]_{u_{_{X_i}}\otimes_RDT} & \mr{Hom}_R(T,A_{X_{i-1}})\ar[d] \\
                A_{X_i}\otimes_RDT\ar[r]          & \mr{Cok}(l(u_{_X}))_i.
        }$ }

\vskip 10pt

Moreover, for  $Y\in \mc{RC}(R)$ and  $h\in\mr{Hom}_{\mc{RC}(R)}(X,Y)$, by applying the homomorphism $l$ to the left square in the commutation diagram in \ref{resolution}, we obtain the following commutative diagram in $\mc{RC}(S)$, for some $h_{\mr{Cok}}$.

\vskip 5pt
\centerline{$\xymatrix{
0\ar[r] & X\hat{\otimes}DT\ar@{->}[rr]^{l(u_{_X})}\ar[d]_{h\hat{\otimes}DT} & & \mr{L}_T(A_X)\ar@{->}[rr]^{\pi_{_{l_X}}}\ar[d]_{\mr{L}_T(h_A)} & & \mr{Cok}(l(u_{_X}))\ar[r]\ar[d]_{h_{\mr{Cok}}} & 0 \\
0\ar[r] & Y\hat{\otimes}DT\ar@{->}[rr]^{l(u_{_Y})} & & \mr{L}_T(A_Y)\ar@{->}[rr]^{\pi_{_{l_Y}}}  & & \mr{Cok}(l(u_{_Y}))\ar[r] & 0 \\
}$}

%
\bskip\
\Exesub\label{functor}\label{Sfct1}
{\it The assignment $\mb{S}_T: \mc{RC}(R)\to \underline{\mc{RC}}(S)$ by setting $X\longmapsto \mr{Cok}(l(u_{_X}))$ is a functor.
}
%
\zskip\

By \ref{coker}, it is sufficient to prove that $\mb{S}_T(h):=h_{\mr{Cok}}=0$ in $\underline{\mc{RC}}(S)$ provided $h=0$. We divide the proof into two steps.

{\it Step 1:} Consider each piece in the commutative diagram in \ref{resolution}. If $h=\{h_i\}=0$, then $h_i=0$ for each $i$. Thus, we have that $u_{_{X_i}}h_{A_i}=0$ and consequently, $h_{A_i}=\pi_{{_{X_i}}}g_i$ for some $g_i: B_{X_i}\to A_{Y_i}$. Since $B_{X_i}\in \mc{B}\subseteq\mc{X}_T$ for each $i$, there are  exact sequences $0\to B_{X_i}\stackrel{b_i}{\longrightarrow}  T_{B_{X_i}}\to B_i'\to 0$ with $T_{B_{X_i}}\in \mr{add}_RT$ and $B_i'\in \mc{B}\subseteq\mr{KerExt}_R^1(-,\mc{A})$. It follows that there exists $t_i\in \mr{Hom}_R(T_{B_{X_i}}, A_{Y_i})$ such that $g_i=b_it_i$. Altogether we obtain the following commutative diagram

\centerline{$\xymatrix{
A_{X_i}\ar[r]^{\pi_{{_{X_i}}}}\ar[d]_{h_{A_i}} &B_{X_i}\ar[d]^{b_i}\ar[dl]^{g_i} \\
A_{Y_i}  & T_{B_{X_i}}\ar[l]^{t_i}.
}$}

\vskip 10pt

\noindent This induces the following commutative diagram in $\mc{RC}^{\mr{tr}}(R)$, where  $T_{B_X}=\{T_{B_{X_i}}\}$.

\vskip 10pt
\centerline{$\xymatrix{
A_X\ar[r]^{\pi_{{_X}}}\ar[d]_{h_A} &B_X\ar[d]^b\ar[dl]^g \\
A_Y  & T_{B_X}\ar[l]^{t}.
}$}

\vskip 10pt

Set $k:= \pi_{{_X}}b$. Then $\mr{L}_T(h_A)=\mr{L}_T(kt)=\mr{L}_T(k)\mr{L}_T(t)$.
\vskip 10pt

{\it Step 2}: Consider the commutative diagram in \ref{coker}.
 Since $u_{_{X_i}}\pi_{{_{X_i}}}=0$, we see that $l(u_{_X})\mr{L}_T(k)=0$ by the definitions. Hence there is some $\theta\in\mr{Hom}_{\mc{RC}(S)}(\mr{Coker}(l(u_{_X})),\mr{L}_T(A_Y))$ such that $\mr{L}_T(k)=\pi_{l_X}\theta$.
Consequently, we have that $\mr{L}_T(h_A)=\mr{L}_T(k)L(t)=\pi_{l_X}\theta \mr{L}_T(t)$. Now we obtain that $\pi_{l_X}h_{\mr{Coker}} =\mr{L}_T(h_{A})\pi_{l_Y}=\pi_{l_X}\theta \mr{L}_T(t)\pi_{{l_Y}}$. Since $\pi_{l_X}$ is epic, we get that $h_{\mr{Cok}} =\mr{L}_T(h_{A})\pi_{l_Y}=\theta \mr{L}_T(t)\pi_{l_Y}$. That is, we have the following commutative diagram.

\vskip 10pt
\centerline{$\xymatrix{
\mr{Cok}({l(u_{_X})})\ar[r]^{\ \ \theta}\ar[d]_{h_{\mr{Cok}}} & \mr{L}_T(T_{B_{X}})\ar[d]^{\mr{L}_T(t)}  \\
\mr{Cok}({l(u_{_Y})})  & L(A_Y)\ar[l]^{\ \ \ \pi_{l_Y}}
}$}

\vskip 10pt
\noindent Note that $\mr{L}_T(T_{B_X})$ is a projective-injective object in $\mc{RC}(S)$, so $h_{\mr{Cok}}=0$ in $\underline{\mc{RC}}(S)$.

%

%
\bskip\
\Exesub\label{Sfct2}
{\it The functor $\mb{S}_T: \underline{\mc{RC}}(R)\to\underline{\mc{RC}}(S)$
}
\zskip\

We will show that the functor $\mb{S}_T$ factors through $\underline{\mc{RC}}(R)$.

 To see this, it is enough to show that $\mb{S}_T(X)$ is a projective object in $\mc{RC}(S)$ whenever $X$ is a projective object in $\mc{RC}(R)$.

 W.l.o.g., we assume that $X=\{X_i\}$ is an indecomposable projective object in $\mc{RC}(R)$.  Thus, 
  $X$ has the form


\vskip 5pt

\setlength{\arraycolsep}{3pt}
{\footnotesize
 \centerline{$\begin{array}{ccccccccccc}\cdots & \rightsquigarrow &  0 & \rightsquigarrow &  \mr{Hom}_R(DR,I) & \stackrel{1}{\rightsquigarrow} &  I & \rightsquigarrow  & 0 & \rightsquigarrow & \cdots\end{array}$}
}
\vskip 5pt

  \noindent where $I$ is on the ($k$-1)-th position, for some $k$.

 Note that $X_k=\mr{Hom}_R(DR,I)\in\mr{add}_RR\subseteq \mc{B}$ and $X_{k-1}=I\in\mr{add}_R(DR)\subseteq   \mc{A}$, so we can choose $A_X$ as  the form

\vskip 5pt
{\footnotesize
 \centerline{$\begin{array}{ccccccccccc}\cdots & \rightsquigarrow &  0 & \rightsquigarrow &  T_k & {\rightsquigarrow} &  I & \rightsquigarrow  & 0 & \rightsquigarrow & \cdots\end{array}$}
}

  \noindent

\noindent where $A_{X_k}=T_k\in\mr{add}_RT$. And we have that the homomorphism $u_{_X}:X\to A_X$ is of the form

\vskip 5pt

{\scriptsize
\noindent\centerline{$\begin{array}{ccccccccccccccc} X:& & \cdots & \rightsquigarrow &  0 & \rightsquigarrow & 0 & \rightsquigarrow &  \mr{Hom}_R(DR,I) & {\rightsquigarrow} &  I & \rightsquigarrow  & 0 & \rightsquigarrow & \cdots \\ \\
\ \ \ \downarrow u_{_X}& & & & & &   & & \ \ \ \ \downarrow u_k && \ \ \ \downarrow 1&&&&\\ \\
A_X:& & \cdots & \rightsquigarrow &  0 & \rightsquigarrow & 0 & \rightsquigarrow &  T_k & {\rightsquigarrow} &  I & \rightsquigarrow  & 0 & \rightsquigarrow & \cdots
\end{array}$
}
}

Then, from the structure of $l(u_{_X})$, we can see that $l(u_{_X})$ is of the form

\vskip 10pt

{\scriptsize
\noindent\centerline{$\begin{array}{ccccccccccccccc} X\hat{\otimes}DT:& & \cdots & \rightsquigarrow &  0 & \rightsquigarrow & 0 & \rightsquigarrow &  \mr{Hom}_R(DR,I)\otimes_RDT & {\rightsquigarrow} &  I\otimes_RDT & \rightsquigarrow  & 0 & \rightsquigarrow & \cdots \\ \\
\ \ \ \downarrow l(u_{_X})& & & & & & \downarrow & & \ \ \ \ \downarrow (-\theta_{l_k},u_k\otimes_DT) && \ \ \ \downarrow 1&&&&\\
 \mr{L}_T(A_X):& & \cdots & \rightsquigarrow &  0 & \rightsquigarrow &\mr{Hom}_R(T,T_k) & \stackrel{(0, \delta_{L_{k+1}})}{\rightsquigarrow} &  \mr{Hom}_R(T,I)\oplus T_k\otimes_RDT & {\rightsquigarrow} &  I\otimes_RDT & \rightsquigarrow  & 0 & \rightsquigarrow & \cdots
\end{array}$
}
}

\noindent where $\theta_{l_k}$ is defined as in  \ref{lpro} and $\delta_{L_{k+1}}$ is defined as in \ref{FcL} repectively. One checks that both homomorphisms $\theta_{l_k}$ and $\delta_{L_{k+1}}$ are in fact isomorphisms. So we obtain that  $S_T(X)=\mr{Coker}(l(u_{_X}))$ is of the form

{\footnotesize
\noindent\centerline{$\begin{array}{ccccccccccccccc}
 & & \cdots & \rightsquigarrow &  0 & \rightsquigarrow &\mr{Hom}_R(T,T_k) & \stackrel{\delta'_{k+1}}{\rightsquigarrow} &  T_k\otimes_RDT & {\rightsquigarrow} &  0& \rightsquigarrow  & 0 & \rightsquigarrow & \cdots
\end{array}$
}}

\vskip 5pt

\noindent %
where $\delta'_{k+1}$ is the induced isomorphism: $\mr{Hom}_R(T,T_k)\otimes_SDS\to T_k\otimes_RDT$. Since $T_k\otimes_RDT\in\mr{add}_S(T\otimes_RDT)=\mr{add}_S(DS)$, we see that $T_k\otimes_RDT$ is an injective $S$-module and that $\mb{S}_T(X)$ is a projective object in $\mc{RC}(S)$.
%

It follows that the functor $\mb{S}_T$ factors through $\underline{\mc{RC}}(R)$. We still denote by $\mb{S}_T$ the induced functor from $\underline{\mc{RC}}(R)$ to $\underline{\mc{RC}}(S)$.


\zskip\

\Exe\label{func2}\label{Q}
{\bf From $\underline{\mc{RC}}$ to $\underline{\mc{RC}}(R)$: The functor $\mb{Q}_{DT}$}
\zskip\
The functor $\mb{Q}_T$ is indeed defined in a way dual to the construction of $\mb{S}_T$.


\Exesub\label{0-w}\label{rFct}
{\it The functor $\mr{R}_{DT}: \mc{RC}^{\mr{tr}}(S)\to \mc{RC}(R)$ 
}

Dually to \ref{FcL}, for any $X=\{X_i\}\in\mc{RC}^{\mr{tr}}(S)$, we define $\mr{R}_{DT}(X)\in\mc{RC}(R)$ as follows.

\vskip 5pt
$(r1)$ the underlying module is $\mr{R}_{DT}(X)_i=\mr{Hom}_S(DT,X_i)\oplus X_{i+1}\otimes_S T$ and,

$(r2)$ the structure map $\delta^{\mr{H}}_i(\mr{R}_{DT})$ is given by $\jz{0}{\delta_{R_i}}{0}{0}$, where
 $\delta_{R_i}$ is the composition

  \hskip 40pt $ \mr{Hom}_S(DT,X_{i})\stackrel{\mr{Hom}_S(DT,\eta^T_{X_{i}})}{\longrightarrow} \mr{Hom}_S(DT,\mr{Hom}_R(T, X_{i}\otimes_S T))$

 \hskip 190pt $\simeq
 \mr{Hom}_R(DT\otimes_S T, X_{i}\otimes_S T)$

  \hskip 190pt $\simeq
 \mr{Hom}_R(DR, X_{i}\otimes_S T)$.

It is easy to see that $\mr{R}_{DT}$ is a functor.

\vskip 10pt

\noindent{\bf Remark}  (1) If $X\in\mc{RC}^{\mr{tr}}(\mr{add}_SDT)$, i.e.,  $X=\{X_i\}$ with each  $X_i\in\mr{add}_SDT$, then $\mr{Hom}_S(DT,X_{i})\in \mr{add}_RR$ and $\delta_{R_i}$ defined above is an isomorphism for each
$i$. It follows that $\mr{R}_{DT}(X)$ is a projective object in $\mc{RC}(R)$ in the case.

(2) In particular, in case $T=S$, we obtain the functor $\mr{R}_{DS}: \mc{RC}^{\mr{tr}}(S)\to\mc{RC}(S)$ which specially send objects in $\mc{RC}^{\mr{tr}}(\mr{add}_SDS)$ to  a projective object in $\mc{RC}(S)$ .

%

\zskip\

\Exesub\label{r-r2}\label{hom}
{\it The functor $\mr{\hat{H}om}(DT,-): \mc{RC}(S)\to\mc{RC}(R)$} 

\sskip\

Let $Y=(Y_i,\delta^{\mr{H}}_i(Y))\in \mc{RC}(S)$. We definie $\mr{\hat{H}om}_S(DT, Y)\in \mc{RC}(R)$ by setting

\sskip\

$(h1)$ the underlying module is $\mr{\hat{H}om}_S(DT, Y)_i=\mr{Hom}_S(DT, Y_i)$ and,

$(h2)$ the structure map $\delta^{\mr{H}}_i(\mr{\hat{H}om}_S(DT, Y))$ is given by the composition:

\sskip\

\hskip 10pt{\footnotesize $\mr{Hom}_S(DT, \delta^{\mr{H}}_i(Y)) \stackrel{\mr{Hom}_S(DT, -y_i)}{\longrightarrow}\mr{Hom}_S(DT, \mr{Hom}_S(DS, Y_{i-1})) \simeq \mr{Hom}_R(DR, \mr{Hom}_S(DT, Y_{i-1})))$.}

\sskip\

Then from the functor property of $\mr{Hom}_S(DT, -)$, one can see that $\mr{\hat{H}om}(DT, -)$ is a functor from $\mc{RC}(S)$ to $\mc{RC}(R)$.

\zskip\

\Exesub\label{wpro}\label{rhom}
{\it The homomorphism $r^X_Y:\mr{Hom}_{\mc{RC}^{\mr{tr}}(S)}(X,Y)\to \mr{Hom}_{\mc{RC}(R)}(\mr{R}_{DT}(X),\mr{\hat{H}om}(DT,Y))$
}
%

Dually to the homomorphism $l^X_Y$, for any $X\in\mc{RC}^{\mr{tr}}(S)$ and $Y\in\mc{RC}(S)$, we have a canonical homomorphism

\vskip 5pt
\centerline{$r^X_Y:\mr{Hom}_{\mc{RC}^{\mr{tr}}(S)}(X,Y)\to \mr{Hom}_{\mc{RC}(R)}(\mr{R}_{DT}(X),\mr{\hat{H}om}(DT,Y))$}
\vskip 5pt

\noindent which is functional in both variables, defined by

\vskip 5pt
\centerline{$r^X_Y: u=\{u_i\}$ \hskip 10pt $\longmapsto$ \hskip 10pt $f=\{f_i\}$ with $f_i=\ljz{\mr{Hom}_S(DT,u_i)}{-\zeta_{r_i}}$,}
\vskip 5pt
%
%
%
\noindent where
$\zeta_{r_i}: X_{i+1}\otimes_ST\to \mr{Hom}_S(DT,Y_i)$ equals to $(u_{i+1}\otimes_ST) \circ (\delta^\mr{H}_{i+1}(Y)\otimes_ST) \circ \epsilon_{\mr{Hom}_S(DT,Y_i)}^T$, i.e., the composition

\vskip 10pt

\hskip 10pt {\footnotesize{$X_{i+1}\otimes_ST \stackrel{u_{i+1}\otimes_ST}{\longrightarrow} Y_{i+1}\otimes_ST \stackrel{\delta^Y_{i+1}\otimes_ST}{\longrightarrow} \mr{Hom}_S(DS,Y_i)\otimes_S T$

\vskip 10pt

\hskip 175pt $\simeq \mr{Hom}_S(T\otimes_RDT, Y_i)\otimes_ST$

\hskip 175pt $\simeq \mr{Hom}_R(T,\mr{Hom}_S(DT,Y_i))\otimes_ST \stackrel{\epsilon_{\mr{Hom}_S(DT,Y_i)}^T}{\longrightarrow}
\mr{Hom}_S(DT,Y_i)$
}}

\vskip 10pt
\noindent{\bf Remark} Using the fact ${_SDS_S}\simeq {_ST\otimes_RDT_S}$ and the adjoint isomorphism \\
\centerline{$\mb{\Gamma}^{DT}: \mr{Hom}_S(X_{i+1}\otimes_SDS,Y_i)\simeq\mr{Hom}_S(X_{i+1}\otimes_ST,\mr{Hom}_S(DT,Y_i))$,}\\ one can easily check that $\zeta_{r_i}$ is the image of the natural homomorphism $(u_{i+1}\otimes_SDS)\circ \delta^{\otimes}_{i+1}(Y)$ under $\Gamma^{DT}$, i.e., $\zeta_{r_i}=\mb{\Gamma}^{DT}((u_{i+1}\otimes_SDS)\circ \delta^{\otimes}_{i+1}(Y))$.

\Exesub\label{epi}
{\it An epimorphism $v_Y: G_Y\to Y$ with $G_Y\in\mc{RC}^{\mr{tr}}(\mc{G})$}, for $Y\in\mc{RC}(S)$.

%
%
\zskip\

%
%
%

Since $(\mc{G},\mc{K})$ is a complete hereditary cotorsion pair in $\mr{mod}S$. It follows that, for any $Y=\{Y_i\}\in\mc{RC}(S)$,  there is an exact sequence $0\to K_Y\stackrel{k_{v_{_Y}}}{\longrightarrow} G_Y\stackrel{v_{_Y}}{\longrightarrow} Y\to 0$ with $K_Y\in\mc{RC}^{\mr{tr}}(\mc{K})$ and $G_Y\in\mc{RC}^{\mr{tr}}(\mc{G})$.

Moreover, for any $h\in\mr{Hom}_{\mc{RC}(S)}(X,Y)$, there is an induced commutative diagram as follows, since $\mr{Ext}_S^1(\mc{G},\mc{K})=0$.

\vskip 10pt
  \centerline{  $\xymatrix{
 0\ar[r] & K_{X}\ar[r]^{k_{v_{_X}}}\ar[d]_{h_K} & G_X\ar[r]^{v_{{_X}}}\ar[d]^{h_G} & X\ar[r]\ar[d]^{h}&0 \\
 0\ar[r] & K_{Y}\ar[r]^{k_{v_{_Y}}}         & G_Y\ar[r]^{v_{{_Y}}}             & Y\ar[r] & 0 \\
        }$ }
\vskip 10pt

\bskip\

\Exesub\label{ker}
{\it The kernel $\mr{Ker}(r(v_Y))$
}
\zskip\

Applying the functor $\mr{Hom}_S(DT,-)$ to the bottom exact sequence in the above diagram, we obtain an induced exact sequence

\vskip 5pt
\centerline{$0\to \mr{Hom}_S(DT,K_Y)\to \mr{Hom}_S(DT,G_Y)\to \mr{Hom}_S(DT,Y)\to 0$,}
\vskip 5pt

\noindent since $K_{Y_i}\in\mc{K}\subseteq \mr{KerExt}_S^1(DT,-)$ for each $i$. Thus, after applying the homomorphism $r$ in \ref{rhom} to the homomorphism $v_X$ in \ref{epi}, we obtain the following exact sequence in  $\mc{RC}(R)$.

 \vskip 5pt
\centerline{$0 \to \mr{Ker}(r(v_Y))\stackrel{\lambda_{r_Y}}{\longrightarrow}  \mr{R}_{DT}(G_Y)\stackrel{r(v_Y)}{\longrightarrow} \mr{\hat{H}om}(DT,Y)\to 0$}
\vskip 5pt

Moreover, for any  $h\in\mr{Hom}_{\mc{RC}(S)}(X,Y)$, by applying the homomorphism $r$ to the right part of the commutative diagram in \ref{epi}, we obtain the following commutative diagram in $\mc{RC}(R)$, for some $h_{\mr{Ker}}$.

\vskip 5pt
\centerline{$\xymatrix{
0\ar[r] & \mr{Ker}(r(v_X))\ar@{->}[rr]^{\lambda_{r_X}}\ar[d]_{h_{\mr{Ker}}} & & \mr{R}_{DT}(G_{X})\ar@{->}[rr]^{r_{v_X}}\ar[d]_{\mr{R}_{DT}(h_{G})} & & \mr{\hat{H}om}(DT,X)\ar[r]\ar[d]_{\mr{\hat{H}om}(DT,h)} & 0 \\
0\ar[r] & \mr{Ker}(r(v_Y))\ar@{->}[rr]^{\lambda_{r_Y}} & & \mr{R}_{DT}(G_{Y})\ar@{->}[rr]^{r_{v_Y}}  & & \mr{\hat{H}om}(DT,Y)\ar[r] & 0 \\
}$}

%
\bskip\
\Exesub\label{Qfct1}
{\it The assignment $\mb{Q}_{DT}: \mc{RC}(S)\to \underline{\mc{RC}}(R)$ by setting $Y\longmapsto \mr{Ker}(r(v_Y))$ is a functor.
}
%
\zskip\

By \ref{ker}, it is sufficient to prove that $\mb{Q}_{DT}(h):=h_{\mr{Ker}}=0$ in $\underline{\mc{RC}}(R)$ provided $h=0$. This is also divided into two steps.

{\it Step 1:} Consider each piece in the commutative diagram in \ref{epi}. If $h=\{h_i\}=0$, then $h_i=0$ for each $i$. Thus, we have that $h_{G_i}v_{Y_i}=0$ and consequently, $h_{G_i}=g_ik_{v_{Y_i}}$ for some $g_i: G_{X_i}\to K_{Y_i}$. Since $K_{Y_i}\in \mc{K}\subseteq {_{DT}\mc{X}}$ for all $i$, there are  exact sequences $0\to K'_{Y_i}\to  DT_{K_{Y_i}}\stackrel{b_i}{\longrightarrow} K_{Y_i}\to 0$ with $DT_{K_{Y_i}}\in \mr{add}_SDT$ and $K'_{Y_i}\in \mc{K}\subseteq\mr{KerExt}_S^1(\mc{G},-)$. It follows that there exists $t_i\in \mr{Hom}_R(G_{Y_i}, DT_{K_{X_i}})$ such that $g_i=t_i b_i$. Altogether we obtain the following commutative diagram.

\centerline{$\xymatrix{
G_{X_i}\ar[r]^{t_i}\ar[dr]^{g_i}\ar[d]_{h_{G_i}} &DT_{K_{Y_i}}\ar[d]^{b_i} \\
G_{Y_i}  & K_{Y_i}\ar[l]^{k_{v_{_{Y_i}}}}
}$}
\vskip 5pt

\noindent It follows that there is a commutative diagram in $\mc{RC}^{\mr{tr}}(S)$, where $DT_{K_Y}:=\{DT_{K_{Y_i}}\}$,

\centerline{$\xymatrix{
G_{X}\ar[r]^{t}\ar[dr]^{g}\ar[d]_{h_{G}} &DT_{K_{Y}}\ar[d]^{b} \\
G_{Y}  & K_{Y}\ar[l]^{k_{v_{_{Y}}}}.
}$}

\vskip 5pt

Set $\beta:=bk_{v_{Y}}$. Then $\mr{R}_{DT}(h_G)=\mr{R}_{DT}(t\beta)=\mr{R}_{DT}(t)\mr{R}_{DT}(\beta)$.

\vskip 5pt
{\it Step 2}: Consider the commutative diagram in \ref{ker}.
 Since $k_{v_{Y_i}}v_{Y_i}=0$, we see that $\mr{R}_{DT}(\beta)r(v_Y)=0$ by the definitions. Hence there is some $\theta\in\mr{Hom}_{\mc{RC}(R)}(\mr{R}_{DT}(G_X),\mr{Ker}(r(v_Y)))$ such that $\mr{R}_{DT}(\beta)=\theta\lambda_{r_Y}$.
Consequently, we have that $\mr{R}_{DT}(h_G)=\mr{R}_{DT}(t)\mr{R}_{DT}(\beta)=\mr{R}_{DT}(t)\theta\lambda_{r_Y}$. Now we obtain that $h_{\mr{Ker}} \lambda_{r_Y} =\lambda_{r_X}\mr{R}_{DT}(h_{G})=\lambda_{r_X}\mr{R}_{DT}(t)\theta\lambda_{r_Y}$. Since $\lambda_{r_Y}$ is monomorphic, we get that $h_{\mr{Ker}} =\lambda_{r_X}\mr{R}_{DT}(t)\theta$. That is, we have the following commutative diagram.

\centerline{$\xymatrix{
\mr{Ker}({r(v_X)})\ar[r]^{\lambda_{r_X}}\ar[d]_{h_{\mr{Ker}}} &\mr{R}_{DT}(G_{Y})\ar[d]^{\mr{R}_{DT}(t)}  \\
\mr{Ker}({r(v_Y)})  & \mr{R}_{DT}(DT_{K_{Y}})\ar[l]^{\theta}
}$}

\vskip 5pt

\noindent Note that $\mr{R}_{DT}(DT_{K_Y})$ is a projective-injective object in $\mc{RC}(R)$, so $h_{\mr{Ker}}=0$ in $\underline{\mc{RC}}(R)$.

%

%
\bskip\
\Exesub\label{Qfct2}
{\it The functor $\mb{Q}_{DT}: \underline{\mc{RC}}(S)\to\underline{\mc{RC}}(R)$
}
\zskip\

We will show that the functor $\mb{Q}_{DT}$ factors through $\underline{\mc{RC}}(S)$.

 To see this, it is enough to show that $\mb{Q}_{DT}(X)$ is a projective object in $\mc{RC}(R)$,  whenever $X$ is a projective object in $\mc{RC}(S)$.

 W.l.o.g., we assume that $X=\{X_i\}$ is an indecomposable projective object in $\mc{RC}(S)$.  Thus, we have that
  $X$ has the form

\vskip 10pt
{\footnotesize
 \centerline{$\begin{array}{ccccccccccc}\cdots & \rightsquigarrow &  0 & \rightsquigarrow &  P & \stackrel{1}{\rightsquigarrow} &  P\otimes_SDS & \rightsquigarrow  & 0 & \rightsquigarrow & \cdots\end{array}$}
}

  \noindent where $P$ is on the ($k$+1)-th position, for some $k$.

 Note that $X_{k+1}=P\in\mr{add}_RR\subseteq \mc{G}$ and $X_{k}=P\otimes_SDS\in\mr{add}_R(DR)\subseteq \mc{K}$, so we can choose $G_X$ as  the form

\vskip 10pt
{\footnotesize
 \centerline{$\begin{array}{ccccccccccc}\cdots & \rightsquigarrow &  0 & \rightsquigarrow & P & {\rightsquigarrow} &   DT_k & \rightsquigarrow  & 0 & \rightsquigarrow & \cdots\end{array}$}
}

  \noindent

\noindent where $G_{X_k}=DT_k\in\mr{add}_RDT$. And we have that the homomorphism $v_X: G_X \to X$ is of the form

\vskip 10pt

{\scriptsize
\noindent\centerline{$\begin{array}{ccccccccccccccc} G_X:& & \cdots & \rightsquigarrow &  0 & \rightsquigarrow & 0 & \rightsquigarrow &  P & {\rightsquigarrow} &  DT_k& \rightsquigarrow  & 0 & \rightsquigarrow & \cdots \\ \\
\ \ \ \downarrow v_X& & & & & &   & & \ \ \ \ \downarrow 1 && \ \ \ \downarrow v_k&&&&\\ \\
X:& & \cdots & \rightsquigarrow &  0 & \rightsquigarrow & 0 & \rightsquigarrow &  P & {\rightsquigarrow} &  P\otimes_SDS & \rightsquigarrow  & 0 & \rightsquigarrow & \cdots
\end{array}$
}
}

\vskip 5pt

Then, from the structure of $r(v_X)$ in \ref{rhom}, we can see that $r(v_X)$ is of the form

\vskip 10pt

{\scriptsize
\noindent\centerline{$\begin{array}{ccccccccccccccc} \mr{R}_{DT}(X):& & \cdots & \rightsquigarrow &  0 & \rightsquigarrow & \mr{Hom}_S(DT,P) & \rightsquigarrow &  \mr{Hom}_S(DT,DT_k)\oplus P\otimes_ST  & \stackrel{(0, \delta_{R_{k}})^{\mr{T}}}{\rightsquigarrow} &  DT_k\otimes_ST & \rightsquigarrow  & 0 & \rightsquigarrow & \cdots \\ \\
\ \ \ \downarrow r(v_X)& & & & & & \downarrow 1 & & \ \ \ \ \downarrow (\mr{Hom}_S(DT,v_k),\zeta_{r_k})^{\mr{T}} &&  \downarrow &&&&\\ \\
 \mr{\hat{H}om}(DT,X):& & \cdots & \rightsquigarrow &  0 & \rightsquigarrow &\mr{Hom}_S(DT,P) & \rightsquigarrow &  \mr{Hom}_S(DT,P\otimes_SDS)& {\rightsquigarrow} &  0 & \rightsquigarrow  & 0 & \rightsquigarrow & \cdots
\end{array}$
}
}

\noindent where $\zeta_{r_k}$ is defined as in  \ref{rhom} and $\delta_{R_{k}}$ is defined as in \ref{rFct} repectively. One checks that both homomorphisms $\zeta_{r_k}$ and $\delta_{R_{k}}$ are in fact isomorphisms. So we obtain that  $\mb{Q}_{DT}(X)=\mr{Ker}(r(v_X))$ is of the form

{\footnotesize
\noindent\centerline{$\begin{array}{ccccccccccccccc}
 & & \cdots & \rightsquigarrow &  0& \rightsquigarrow & 0 & \rightsquigarrow &\mr{Hom}_S(DT,DT_k) & \stackrel{\delta'_{k}}{\rightsquigarrow} &  DT_k\otimes_ST & {\rightsquigarrow} &  0& \rightsquigarrow & \cdots
\end{array}$
}}

\vskip 5pt

\noindent %
where $\delta'_{k}$ is an induced isomorphism: $\mr{Hom}_S(DT,DT_k) \to \mr{Hom}_S(DS,DT_k\otimes_ST)$. Since $\mr{Hom}_S(DT,DT_k)\in\mr{add}_R(\mr{Hom}_S(DT,DT))=\mr{add}_RR$, we see that $\mr{Hom}_S(DT,DT_k)$ is projective and that $\mb{Q}_{DT}(X)$ is a projective object in $\mc{RC}(R)$.
%

It follows that the functor $\mb{Q}_{DT}$ factors through $\underline{\mc{RC}}(S)$. We still denote by $\mb{Q}_{DT}$ the induced functor from $\underline{\mc{RC}}(S)$ to $\underline{\mc{RC}}(R)$.


\zskip\

\Exe\label{QS=1}
{\bf The isomorphism  $\mb{Q}_{DT}\mb{S}_T\simeq 1_{\underline{\mc{RC}}(R)}$}
%
%
%
\Exesub\label{QSX}
{\it Computing the composition $\mb{Q}_{DT}\mb{S}_T$ 
}
\zskip\

Take any $X=\{X_i,\delta^{\otimes}_i\}\in\mc{RC}(R)$.
From the chosen exact sequence $0\to X\stackrel{u_{_X}}{\longrightarrow} A_X\stackrel{\pi_{u_{_X}}}{\longrightarrow} B_X\to 0$ in $\mc{RC}^{\mr{tr}}(R)$ with $A_X=\{A_{X_i}\}\in\mc{RC}^{\mr{tr}}(\mc{A})$ and $B_X=\{B_{X_i}\}\in\mc{RC}^{\mr{tr}}(\mc{B})$, as in \ref{resolution}, we obtain an exact sequence

\vskip 5pt
\centerline{$0\to X\hat{\otimes} DT\stackrel{l(u_{_X})}{\longrightarrow} \mr{L}_T(A_X)\stackrel{s}{\longrightarrow} \mb{S}_T(X)\to 0$}
\vskip 5pt

\noindent by the construction of the functor $\mb{S}_T$ in \ref{coker}. 
 Note that, for each $i$, $\mb{S}_T(X)_i$ is given by the pushout diagram

  \centerline{  $\xymatrix{
             {X_i}\otimes_RDT\ar[r]^{\theta_{l_i}\  \  }\ar[d]_{u_{_{X_i}}\otimes_RDT} & \mr{Hom}_R(T,A_{X_{i-1}})\ar[d]_{s^1_i} \\
                A_{X_i}\otimes_RDT\ar[r]^{s^2_i}          & \mb{S}_T(X)_i.
        }$ }

\vskip 10pt

Now we take a projective $R$-module $P_{A_{X_i}}$ such that $P_{A_{X_i}}\stackrel{p_i}{\longrightarrow} A_{X_i}\to 0$ is exact.  Then we have a pullback diagram

\vskip 10pt
  \centerline{  $\xymatrix{
 0\ar[r] & \overline{X}_i\ar[r]^{\overline{u}_{_{X_i}}}\ar[d]_{q_i} & P_{A_{X_i}}\ar[r]\ar[d]^{p_i} & B_{X_i}\ar[r]\ar@{= }[d]&0 \\
 0\ar[r] & X_i\ar[r]^{u_{_{X_i}}}         & A_{X_i}\ar[r]           & B_{X_i}\ar[r] & 0. \\
        }$ }
\vskip 10pt

Since $\mathcal{B}$ is  closed under kernels of epimorphisms, we see that $\overline{X}_i\in\mathcal{B}$.
By applying the functor $-\otimes_R DT$, the diagram above induces the following commutative diagram with exact rows.

\vskip 5pt
\centerline{$\xymatrix{
0\ar[r] & \overline{X}_i\otimes_R DT\ar@{->}[rr]^{\overline{u}_{_{X_i}}\otimes_R DT}\ar[d]_{q_i\otimes_R DT} & & P_{A_{X_i}}\otimes_R DT\ar@{->}[rr] \ar[d]_{p_i\otimes_R DT} & & B_{X_i}\otimes_R DT\ar[r]\ar@{= }[d]  & 0 \\
0\ar[r] & X_i\otimes_R DT\ar@{->}[rr]^{u_{_{X_i}}\otimes_R DT} & & A_{X_i}\otimes_R DT\ar@{->}[rr]  & & B_{X_i}\otimes_R DT\ar[r] & 0 \\
}$}

\vskip 10pt

Now, one can check that  the following diagram is commutative with exact rows, for each $i$, where the lower row is obtained from the above exact sequence.

%
%

%
\bskip\
 \setlength{\unitlength}{0.09in}
\hskip 60pt
 \begin{picture}(50,7)

 \put(-9,0){\makebox(0,0)[c]{$0$}}
                             \put(-8,0){\vector(1,0){2}}

 \put(-2,0){\makebox(0,0)[c]{$X_i\otimes_R DT$}}
                             \put(8,1){\makebox(0,0)[c]{{\footnotesize $_{(-\theta_{l_i},u_{_{X_i}}\otimes_RDT)}$}}}
                             \put(3,0){\vector(1,0){10}}
                 \put(-2,4){\vector(0,-1){2}}
                       \put(1,3){\makebox(0,0)[c]{$_{q_i\otimes_R DT}$}}
 \put(25,0){\makebox(0,0)[c]{{\footnotesize $\mr{Hom}_R(T,A_{X_{i-1}})\oplus A_{X_i}\otimes_R DT$}}}  
                       \put(38,0){\vector(1,0){7}}
                            \put(42,2){\makebox(0,0)[c]{{\scriptsize ${s_i=\ljz{s^1_i}{s^2_i}}$}}}
\put(25,4){\vector(0,-1){2}}
                       \put(30,3){\makebox(0,0)[c]{{\tiny $_{\jz{1}{0}{0}{p_i\otimes_R DT}}$}}}

\put(50,0){\makebox(0,0)[c]{$\mb{S}_T(X)_i$}}
                            \put(53,0){\vector(1,0){2}}
\put(50,4){\line(0,-1){2}}
\put(50.2,4){\line(0,-1){2}}
 \put(56,0){\makebox(0,0)[c]{$0$}}

 \put(-9,6){\makebox(0,0)[c]{$0$}}
                             \put(-8,6){\vector(1,0){2}}
 \put(-2,6){\makebox(0,0)[c]{$\overline{X}_i\otimes_R DT$}}
                             \put(8,7){\makebox(0,0)[c]{{\footnotesize $_{(-(q_i\otimes_R DT)\circ \theta_{l_i},\overline{u}_{X_i}\otimes_R DT)}$}}}
                             \put(3,6){\vector(1,0){10}}
 \put(25,6){\makebox(0,0)[c]{{\footnotesize $\mr{Hom}_R(T,A_{X_{i-1}})\oplus P_{A_{X_i}}\otimes_R DT$}}}

                             \put(38,6){\vector(1,0){7}}
                              \put(42,8){\makebox(0,0)[c]{{\tiny ${s^P_i=\ljz{s^1_i}{(p_i\otimes_R DT)\circ s^2_i}}$}}}
 \put(50,6){\makebox(0,0)[c]{$\mb{S}_T(X)_i$}}
                             \put(53,6){\vector(1,0){2}}
 \put(56,6){\makebox(0,0)[c]{$0$}}

\end{picture}
\bskip\

Denote $\mf{L}^P_{A_X}:=\{\mr{Hom}_R(T,A_{X_{i-1}})\oplus P_{A_{X_i}}\otimes_R DT\}\in\mc{RC}^{\mr{tr}}(S)$. Note that $\overline{X}_i\otimes_R DT\in \mc{K}$ and $\mr{Hom}_R(T,A_{X_{i-1}})\oplus P_{A_{X_i}}\otimes_R DT\in \mc{G}$, so we have an exact sequence from the first row in the commutative diagram 

\vskip 5pt
\centerline{$0\to \overline{X}\otimes_R DT{\longrightarrow} \mf{L}^P_{A_X}\stackrel{s^P}{\longrightarrow} \mb{S}_T(X)\to 0$}
\vskip 5pt 
 
\noindent with $\overline{X}\otimes_R DT\in\mc{RC}^{\mr{tr}}(\mc{K})$ and $\mf{L}^P_{A_X}\in\mc{RC}^{\mr{tr}}(\mc{G})$, as in \ref{epi}. By applying the homomorphism $r$ in \ref{rhom} to the homomorphism $s^P: \mf{L}^P_{A_X}\to \mb{S}_T(X)$, we have an exact sequence in $\mc{RC}(R)$ by the construction of the functor $\mb{Q}_{DT}$ in \ref{ker}

\vskip 5pt
\centerline{$0\to \mb{Q}_{DT}\mb{S}_T(X)\stackrel{\lambda}{\longrightarrow} \mr{R}_{DT}(\mf{L}^P_{A_X})\stackrel{r(s^P)}{\longrightarrow} \mr{\hat{H}om}(DT, \mb{S}_T(X))\to 0$,}
\vskip 5pt

\noindent where $r(s^P)$ is defined as in \ref{rhom}.


\zskip\

\Exesub\label{X+P}
{\it The object $X\oplus  \mr{L}_R(P^+_{A_X})$ in $\mc{RC}(R)$}
\zskip\

Denote that $P^+_{A_X}:=\{P_{A_{X_{i+1}}}\}$, then $P^+_{A_X}\in\mc{RC}^{\mr{tr}}(\mr{add}_RR)$. Applying the functor  $\mr{L}_R$ in the remark in \ref{FcL}, we obtain that  $\mr{L}_R(P^+_{A_X})$ is a projective object in $\mc{RC}(R)$. Hence, the object $X\oplus  \mr{L}_R(P^+_{A_X})$ is isomorphic to $X$ in $\underline{\mc{RC}}(R)$.

We will prove that $\mb{Q}_{DT}\mb{S}_T(X)\simeq X\oplus  \mr{L}_R(P^+_{A_X})$ naturally. And then, $\mb{Q}_{DT}\mb{S}_T \simeq 1_{\underline{\mc{RC}}(R)}$.

The general strategy is as follows. Firstly, we construct a natural homomorphism $\xi: X\oplus \mr{L}_R(P^+_{A_X})\to \mr{R}_{DT}(\mf{L}^P_{A_X})$. Secondly, we show that $\xi\circ r(s^P)=0$, i.e.,  the composition of $\xi$ and the homomorphism $r(s^P): \mr{R}_T(\mf{L}^P_{A_X})\to \mr{\hat{H}om}(DT, \mb{S}_T(X))$ in the exact sequence above is 0. Thus, we obtain a homomorphism $\phi: X\oplus \mr{L}_R(P^+_{A_X})\to \mb{Q}_{DT}\mb{S}_T(X)$. Finally, we prove that $\phi$ is indeed a  natural isomorphism.

%

\vskip 10pt

\Exesub\label{theta}
{\it   The homomorphism $\xi: X\oplus \mr{L}_R(P^+_{A_X})\to \mr{R}_{DT}(\mf{L}^P_{A_X})$
}

\vskip 10pt

Recall from the construction in \ref{QSX} that

    \hskip 30pt $X=\{X_i\}$ and,

    \hskip 30pt $\mr{L}_R(P^+_{A_X})=\{\mr{Hom}_R(R,P_{A_{X_{i}}})\oplus P_{A_{X_{i+1}}}\otimes_RDR\}=\{P_{A_{X_{i}}}\oplus P_{A_{X_{i+1}}}\otimes_RDR\}$

\noindent and that

\hangafter 1\hangindent 97pt
       \hskip 30pt  $\mr{R}_{DT}(\mf{L}^P_{A_X})=\{\mr{Hom}_S(DT, (\mf{L}^P_{A_X})_i)\oplus  (\mf{L}^P_{A_X})_{i+1}\otimes_S T\}$\\ $=\{\mr{Hom}_S(DT, \mr{Hom}_R(T,A_{X_{i-1}})\oplus P_{A_{X_i}}\otimes_R DT)\oplus  (\mr{Hom}_R(T,A_{X_i})\oplus  P_{A_{X_{i+1}}}\otimes_R DT)\otimes_S T\}$

Let  $\xi=\{\xi_i\}: X\oplus \mr{L}_R(P^+_{A_X})\to \mr{R}_{DT}(\mf{L}^P_{A_X})$ be a homomorphism. We may assume that $\xi_i=(\xi_i^a,\xi_i^b)$, where

       \hskip 30pt  $\xi_i^a: X_i\oplus P_{A_{X_{i}}}\oplus P_{A_{X_{i+1}}}\otimes_RDR\to \mr{Hom}_S(DT, (\mf{L}^P_{A_X})_i)$ and

       \hskip 30pt  $\xi_i^b: X_i\oplus P_{A_{X_{i}}}\oplus P_{A_{X_{i+1}}}\otimes_RDR\to (\mf{L}^P_{A_X})_{i+1}\otimes_ST$.

%

\vskip 5pt
\Exesubsub\label{xi-a}
{\it The homomorphism $\xi_i^a$ in $\mr{mod}R$
}

 \vskip 10pt

We set {\footnotesize $\xi_i^a=\32jz{\xi^a_{11}}{\xi^a_{12}}{\xi^a_{21}}{\xi^a_{22}}{\xi^a_{31}}{\xi^a_{32}}:$}

\vskip 10pt
\centerline{ $X_i\oplus P_{A_{X_{i}}}\oplus P_{A_{X_{i+1}}}\otimes_RDR\to \mr{Hom}_S(DT, \mr{Hom}_R(T,A_{X_{i-1}})\oplus P_{A_{X_i}}\otimes_R DT)$.}

\vskip 10pt
\noindent Using the isomorphism ${_SDS_S}\simeq {_ST\otimes_RDT_S}$ and the adjoint isomorphism

\vskip 10pt
\centerline{$\mb{\Gamma}^{DT}: \mr{Hom}_S(-\otimes_SDT,-)\simeq\mr{Hom}_S(-,\mr{Hom}_S(DT,-))$,}

\vskip 10pt
\noindent we define the components as follows.

\vskip 10pt
\noindent $\bullet$ We set

\vskip 10pt
\centerline{$\xi^a_{11}=\mb{\Gamma}^{DT}(\theta_{l_i}): X_i\to \mr{Hom}_S(DT,\mr{Hom}_R(T,A_{X_{i-1}}))$, }

\vskip 10pt

\noindent where $\theta_{l_i}: X_i\otimes_SDT\to\mr{Hom}_R(T,A_{X_{i-1}})$ is defined in \ref{lhom}.

In the other words, the morphism  $\xi^a_{11}$ is given by the composition : $\eta_{_{X_i}}^{DR}\circ \mr{Hom}_R(DR,\delta_{X_i}) \circ \mr{Hom}_R(DR,\ u_{_{X_{i-1}}})$ and some natural isomorphisms

 \bg{verse}
 \bg{verse}

 $X_i \stackrel{\eta_{_{X_i}}^{DR}}{\longrightarrow} \mr{Hom}_R(DR,X_i\otimes_RDR) \stackrel{\mr{Hom}_R(DR,\delta_{X_i})}{\longrightarrow} \mr{Hom}_R(DR,X_{i-1})$

   $  \stackrel{\mr{Hom}_R(DR,\ u_{_{X_{i-1}}})}{\longrightarrow}  \mr{Hom}_R(DR,A_{X_{i-1}})  \simeq \mr{Hom}_R(DT\otimes_ST, A_{X_{i-1}})$

    \hskip 162pt $ \simeq \mr{Hom}_R(DT,\mr{Hom}_R(T, A_{X_{i-1}}))$.

\end{verse}

\end{verse}

 $\bullet$ The morphism  $\xi^a_{22}=\mb{\Gamma}^{DT}(1_{(P_{A_{X_i}}\otimes_R DT)})=\eta^{DT}_{_{P_{A_{X_i}}}}: P_{A_{X_i}}\to \mr{Hom}_S(DT,P_{A_{X_i}}\otimes_R DT)$.

 $\bullet$ The remained morphisms $\xi^a_{12}, \xi^a_{21},\xi^a_{31},\xi^a_{32}$ are all $0$.

 \vskip 10pt

 So we have that  {\footnotesize $\xi_i^a=\32jz{\xi^a_{11}}{0}{0}{\xi^a_{22}}{0}{0}$}, where $\xi^a_{11}=\mb{\Gamma}^{DT}(\theta_{l_i})$, and $\xi^a_{22}=\mb{\Gamma}^{DT}(1_{(P_{A_{X_i}}\otimes_R DT)})$.

%
%
%
%
%
%

\vskip 10pt

\vskip 5pt
\Exesubsub\label{xi-b}
{\it The homomorphism $\xi_i^b$ in $\mr{mod}R$
}

 \vskip 10pt

We set {\footnotesize $\xi_i^b=\32jz{\xi^b_{11}}{\xi^b_{12}}{\xi^b_{21}}{\xi^b_{22}}{\xi^b_{31}}{\xi^b_{32}}:$}
\vskip 5pt

\centerline{ $X_i\oplus P_{A_{X_{i}}}\oplus P_{A_{X_{i+1}}}\otimes_RDR\to (\mr{Hom}_R(T,A_{X_i})\oplus  P_{A_{X_{i+1}}}\otimes_R DT)\otimes_S T$,}

\noindent where the components are defined naturally as follows.

 \vskip 10pt

\hangafter 1\hangindent 30pt
 $\bullet$  The morphism $\xi^b_{11}=u_{_{X_i}}\circ (\epsilon^T_{A_{X_i}})^{-1}: X_i\to \mr{Hom}_R(T,A_{X_i})\otimes_S T$ (note that $\epsilon^T_{A_{X_i}}$ is an isomorphism since $A_{X_i}\in\mc{A}$), i.e., is given by the composition

      \centerline{$X_i\stackrel{u_{_{X_i}}}{\longrightarrow} A_{X_i}\stackrel{(\epsilon^T_{A_{X_i}})^{-1}}{\longrightarrow} \mr{Hom}_R(T,A_{X_i})\otimes_S T$.}

\hangafter 1\hangindent 30pt
$\bullet$  The morphism $\xi^b_{21}=p_i\circ (\epsilon^T_{A_{X_i}})^{-1}: P_{A_{X_i}}\to \mr{Hom}_R(T,A_{X_i})\otimes_S T$, i.e., is given by the composition

      \centerline{$P_{A_{X_i}}\stackrel{p_i}{\longrightarrow} A_{X_i}\stackrel{(\epsilon^T_{A_{X_i}})^{-1}}{\longrightarrow} \mr{Hom}_R(T,A_{X_i})\otimes_S T$.}

\hangafter 1\hangindent 30pt
 $\bullet$
     The morphism $\xi^b_{32}: P_{A_{X_{i+1}}}\otimes_R DR\to P_{A_{X_{i+1}}}\otimes_R DT\otimes_S T$ is the natural isomorphism given by ${_R(DT\otimes_ST)_R}\simeq {_RR_R}$.

 $\bullet$ The remained morphisms $\xi^b_{12}, \xi^b_{22},\xi^b_{31}$ are all $0$.

 \vskip 10pt

 So we have that  {\footnotesize $\xi_i^b=\32jz{\xi^b_{11}}{0}{\xi^b_{21}}{0}{0}{\xi^b_{32}}$}, where $\xi^b_{11}=u_{_{X_i}}\circ (\epsilon^T_{A_{X_i}})^{-1}$, $\xi^b_{21}=p_i\circ (\epsilon^T_{A_{X_i}})^{-1}$, and $\xi^b_{32}$ is the natural isomorphism.

%
%
%

\vskip 10pt

\Exesubsub\label{xi-hom}
{\it $\xi$ is a homomorphism in $\mc{RC}(R)$
}

 \vskip 10pt

It is not difficult to prove that the above-defined  morphism $\xi$ is in fact a homomorphism in $\mc{RC}(R)$ by the involved definitions.
%
%
%
%
%
%
%
%
%

\vskip 10pt

\Exesub\label{xi-rsp=0}
\hangafter 1\hangindent 40pt
{\it The composition $\xi\circ r(s^P)=0$, and so $\xi$ factors through a homomorphism $\phi: X\oplus \mr{L}_R(P^+_{A_X})\to \mb{Q}_{DT}\mb{S}_T(X)$.
}

 \vskip 10pt

\Exesubsub\label{s-fx}
\hangafter 1\hangindent 40pt
{\it The analysis of the homomorphism $s: \mr{L}_T(A_X)\to \mb{S}_T(X)$ in \ref{QSX}.
}

 \vskip 10pt

Recall from \ref{QSX} that $s=\{s_i\}: \mr{L}_T(A_X)\to \mb{S}_T(X)$ is a homomorphism in $\mc{RC}(S)$ which is the cokernel of  the homomorphism $l(u_{_X})$. Note that $\mr{L}_T(A_X)=\mr{Hom}_R(T,A_{X_{i-1}})\oplus A_{X_i}\otimes DT$, so we write that $s_i=\ljz{s^1_i}{s^2_i}$ as we have done in the last commutative diagram in \ref{QSX}. The fact that $s$ is a  homomorphism in $\mc{RC}(S)$ implies that there is the following commutative diagram, for each $i$.

\bskip\

 \setlength{\unitlength}{0.09in}
 \begin{picture}(50,7)

 \put(17,0){\makebox(0,0)[c]{$_{\mr{Hom}_R(T,A_{X_{i-2}})\oplus A_{X_{i-1}}\otimes_R DT}$}}
                             \put(38,2){\makebox(0,0)[c]{{\tiny $\ljz{s^1_{i-1}}{s^2_{i-1}}$}}}
                             \put(33,0){\vector(1,0){8}}
 \put(47,0){\makebox(0,0)[c]{$_{\mb{S}_T(X)_{i-1}}$}}

                 \put(17,4){\vector(0,-1){2}}
                       \put(22,3){\makebox(0,0)[c]{$_{\delta^{\otimes}_i(\mr{L}_T(A_X))}$}}
                 \put(47,4){\vector(0,-1){2}}
                       \put(51,3){\makebox(0,0)[c]{$_{\delta^{\otimes}_i({\mb{S}_T(X)})}$}}

 \put(17,6){\makebox(0,0)[c]{$_{[\mr{Hom}_R(T,A_{X_{i-1}})\oplus A_{X_i}\otimes_R DT]\otimes_S DS}$}}
                             \put(38,8){\makebox(0,0)[c]{{\tiny $\ljz{s^1_{i}}{s^2_{i}}\otimes_S DS$}}}
                             \put(33,6){\vector(1,0){8}}
 \put(47,6){\makebox(0,0)[c]{$_{\mb{S}_T(X)_i\otimes_S DS}$}}

\end{picture}
\bskip\
%
{ 
\hskip 15pt By the definition of $\delta^{\otimes}_i(\mr{L}_T(A_X))$ (see \ref{FcL}) and the above commutative diagram, we obtain that $(s^2_i\otimes DS)\circ \delta^{\otimes}_i({\mb{S}_T(X)})=0$ and that $(s^1_i\otimes DS)\circ  \delta^{\otimes}_i({\mb{S}_T(X)})=(\epsilon^T_{A_{X_{i-1}}}\otimes DT) \circ s^2_{i-1}$.}

%
%
%

\Exesubsub\label{rsp-fx}
\hangafter 1\hangindent 40pt
{\it The analysis of the homomorphism $r(s^P): \mr{R}_T(\mf{L}^P_{A_X})\to \mr{\hat{H}om}(DT,\mb{S}_T(X))$
}

 \vskip 10pt

 Recall from \ref{QSX} that $s^P=\{s^P_i\}=\{\ljz{s^1_i}{(p_i\otimes_RDT)\circ s^2_i}\}: \mf{L}^P_{A_X}\to \mb{S}_T(X)$. By \ref{rhom}, we know that $\mr{R}_T(\mf{L}^P_{A_X})=\{\mr{Hom}_S(DT,(\mf{L}^P_{A_X})_i)\oplus (\mf{L}^P_{A_X})_{i+1}\otimes_ST\}$ and that $r(s^P)_i=\ljz{\mr{Hom}_S(DT,s^P_i)}{-\zeta_{r_i}}$, where $\zeta_{r_i}=\mb{\Gamma}^{DT}((s^P_{i+1}\otimes_SDS)\circ \delta^{\otimes}_{i+1}(\mb{S}_T(X)))$.

%
%
%
%
%
%
%
%
%
%

For convenience, we set $r(s^P)_i=\ljz{r^1_i}{r^2_i}$, where

\vskip 10pt
\hskip 50pt $r^1_i=\mr{Hom}_S(DT,s^P_i): \mr{Hom}_S(DT,(\mf{L}^P_{A_X})_i)\to \mb{S}_T(X)_i$ and

\vskip 10pt
\hskip 50pt $r^2_i=-\zeta_{r_i}: (\mf{L}^P_{A_X})_{i+1}\otimes_ST\to \mb{S}_T(X)_i$.


\vskip 10pt
\Exesubsub\label{xi-rsp=0check}
\hangafter 1\hangindent 40pt
{\it Checking $\xi\circ r(s^P)=0$
}

\vskip 10pt

To check  $\xi\circ r(s^P)=0$, we need only to check that $\xi_i^ar^1_i+\xi_i^br^2_i=0$ for each $i$, since $\xi_i=(\xi_i^a,\xi_i^b)$ and $r(s^P)_i=\ljz{r^1_i}{r^2_i}$. Note that $r^1_i=\mr{Hom}_S(DT,s^P_i)$ and $r^2_i=-\zeta_{r_i}$, so it is enough to check that $\xi^a_i\circ \mr{Hom}_S(DT,s^P_i)=\xi_i^b\circ \zeta_{r_i}$.

\vskip 10pt
Since {\footnotesize $\xi_i^a=\32jz{\xi^a_{11}}{0}{0}{\xi^a_{22}}{0}{0}$},
      {\footnotesize $\xi_i^b=\32jz{\xi^b_{11}}{0}{\xi^b_{21}}{0}{0}{\xi^b_{32}}$},
      $\zeta_{r_i}=\mb{\Gamma}^{DT}((s^P_{i+1}\otimes_SDS)\circ \delta^{\otimes}_{i+1}(\mb{S}_T(X)))$ and
       $s^P_i=\ljz{s^1_i}{(p_i\otimes_RDT)\circ s^2_i}$, we just check the following.

 \vskip 10pt
 \vskip 10pt
\noindent $(1)$ $\xi^a_{11}\circ\mr{Hom}_S(DT,s^1_i)=\xi^b_{11}\circ \mb{\Gamma}^{DT}((s^1_{i+1}\otimes_SDS)\circ \delta^{\otimes}_{i+1}(\mb{S}_T(X)))$.

 \vskip 10pt
\hskip 20pt
\hangafter 1\hangindent 20pt 
By \ref{xi-a}, we have that

 \vskip 10pt
\hskip 40pt
\hangafter 1\hangindent 45pt  $\xi^a_{11}\circ\mr{Hom}_S(DT,s^1_i)=\mb{\Gamma}^{DT}(\theta_{l_i})\circ\mr{Hom}_S(DT,s^1_i)=\mb{\Gamma}^{DT}(\theta_{l_i}\circ s^1_i)$,
 \vskip 10pt

\hangafter 1\hangindent 20pt
where later equality uses the naturality of $\mb{\Gamma}^{DT}$. On the other hand, by \ref{xi-b} and \ref{s-fx}, we obtain that

 \vskip 10pt
\hskip 40pt
\hangafter 1\hangindent 45pt  $\xi^b_{11}\circ \mb{\Gamma}^{DT}((s^1_{i+1}\otimes_SDS)\circ \delta^{\otimes}_{i+1}(\mb{S}_T(X)))$\\
  $ =(u_{_{X_i}}\circ (\epsilon^T_{A_{X_i}})^{-1})\circ  \mb{\Gamma}^{DT}((\epsilon^T_{A_{X_i}}\otimes_R DT) \circ s^2_i)$\\
  $ =\mb{\Gamma}^{DT}(((u_{_{X_i}}\circ (\epsilon^T_{A_{X_i}})^{-1})\otimes_RDT)\circ (\epsilon^T_{A_{X_i}}\otimes_R DT) \circ s^2_i)$\\
  $ =\mb{\Gamma}^{DT}((u_{_{X_i}}\otimes_RDT)\circ  s^2_i)$.
 \vskip 10pt

\hangafter 1\hangindent 20pt
But $(u_{_{X_i}}\otimes_RDT)\circ  s^2_i= \theta_{l_i}\circ s^1_i$ by the pushout diagram on $\mb{S}_T(X)_i$ in \ref{QSX}. Hence, the equality (1) holds.

 \vskip 10pt
\noindent $(2)$ $\xi^a_{22}\circ\mr{Hom}_S(DT,(p_i\otimes_RDT)\circ s^2_i)=\xi^b_{21}\circ \mb{\Gamma}^{DT}((s^1_{i+1}\otimes_SDS)\circ \delta^{\otimes}_{i+1}(\mb{S}_T(X)))$.

\vskip 10pt
\hskip 20pt
\hangafter 1\hangindent 20pt 
By \ref{xi-a} and  the naturality of $\mb{\Gamma}^{DT}$,

 \vskip 10pt
\hskip 40pt
\hangafter 1\hangindent 45pt  $\xi^a_{22}\circ\mr{Hom}_S(DT,(p_i\otimes_RDT)\circ s^2_i)$\\
$=\mb{\Gamma}^{DT}(1_{P_{A_{X_i}}\otimes_RDT})\circ\mr{Hom}_S(DT,(p_i\otimes_RDT)\circ s^2_i)$\\
$=\mb{\Gamma}^{DT}(1_{P_{A_{X_i}}\otimes_RDT}\circ ((p_i\otimes_RDT)\circ s^2_i))$\\
$=\mb{\Gamma}^{DT}((p_i\otimes_RDT)\circ s^2_i)$.

\vskip 10pt
\hangafter 1\hangindent 20pt
On the other hand, by \ref{xi-b} and \ref{s-fx} and   the naturality of $\mb{\Gamma}^{DT}$,

 \vskip 10pt
\hskip 40pt
\hangafter 1\hangindent 45pt  $\xi^b_{21}\circ \mb{\Gamma}^{DT}((s^1_{i+1}\otimes_SDS)\circ \delta^{\otimes}_{i+1}(\mb{S}_T(X)))$\\
$=(p_i\circ (\epsilon^T_{A_{X_i}})^{-1})\circ  \mb{\Gamma}^{DT}((\epsilon^T_{A_{X_i}}\otimes_R DT) \circ s^2_i)$\\
  $ =\mb{\Gamma}^{DT}(((p_i\circ (\epsilon^T_{A_{X_i}})^{-1})\otimes_RDT)\circ (\epsilon^T_{A_{X_i}}\otimes_R DT) \circ s^2_i)$\\
  $=\mb{\Gamma}^{DT}((p_i\otimes_RDT)\circ  s^2_i)$.

 \vskip 10pt
\hangafter 1\hangindent 20pt
Hence, the equality (2) holds.

 \vskip 10pt
\noindent $(3)$ $0=\xi^b_{32}\circ \mb{\Gamma}^{DT}(((p_{i+1}\otimes_RDT)\circ s^2_{i+1})\otimes_SDS)\circ \delta^{\otimes}_{i+1}(\mb{S}_T(X)))$.

 \vskip 10pt
\hskip 20pt
\hangafter 1\hangindent 20pt
In fact, the equality holds by observing that

 \vskip 10pt
\hskip 40pt
\hangafter 1\hangindent 45pt  $\mb{\Gamma}^{DT}(((p_{i+1}\otimes_RDT)\circ s^2_{i+1})\otimes_SDS)\circ \delta^{\otimes}_{i+1}(\mb{S}_T(X)))$\\
$=\mb{\Gamma}^{DT}((p_{i+1}\otimes_RDT\otimes_SDS) \circ (s^2_{i+1}\otimes_SDS)\circ \delta^{\otimes}_{i+1}(\mb{S}_T(X)))$\\
$=0$,

 \vskip 10pt
\hangafter 1\hangindent 20pt
 since $(s^2_{i+1}\otimes_SDS)\circ \delta^{\otimes}_{i+1}(\mb{S}_T(X)))=0$ by \ref{s-fx}.

 \vskip 10pt
All together, we prove that $\xi\circ r(s^P)=0$ and therefore, $\xi$ factors through $\mb{Q}_{DT}\mb{S}_T(X)=\mr{Ker}(r(s^P))$ by a homomorphism  $\phi: X\oplus \mr{L}_R(P^+_{A_X})\to \mb{Q}_{DT}\mb{S}_T(X)$ in $\mc{RC}(R)$, i.e., $\xi=\phi\circ \lambda$.


\vskip 10pt
\Exesub\label{xi-iso}
\hangafter 1\hangindent 40pt
{\it The induced homomorphism $\phi: X\oplus \mr{L}_R(P^+_{A_X})\to \mb{Q}_{DT}\mb{S}_T(X)$ is an isomorphism
}

\vskip 10pt

We now prove that the induced  homomorphism $\phi: X\oplus \mr{L}_R(P^+_{A_X})\to \mb{Q}_{DT}\mb{S}_T(X)$ is an isomorphism. Clearly it is equivalent to show that $\phi_i: (X\oplus \mr{L}_R(P^+_{A_X}))_i\to \mb{Q}_{DT}\mb{S}_T(X)_i$ is an isomorphism,  for each $i$.

We will show that  there is the following commutative diagram ($\ast$) with exact rows, for each $i$, where $\overline{X}_i\in\mc{B}$ is obtained in \ref{QSX}. Note that $\mr{L}_R(P^+_{A_X})_i=\mr{Hom}_R(T,A_{X_i})\otimes_S T\oplus P_{A_{X_{i+1}}}\otimes_R DT\otimes_S T\simeq A_{X_i}\oplus P_{A_{X_{i+1}}}\otimes_R DR$.

\vskip 10pt
\bskip\

 \setlength{\unitlength}{0.09in}
 \begin{picture}(50,7)

                 \put(-2,3){\makebox(0,0)[c]{$(\ast):$}}

 \put(0,0){\makebox(0,0)[c]{$_0$}}
                             \put(1,0){\vector(1,0){2}}

 \put(10,0){\makebox(0,0)[c]{$_{\mr{Hom}_R(DT, \overline{X}_i\otimes_R DT)}$}}
                             \put(19,1){\makebox(0,0)[c]{$_{a_i}$}}
                             \put(17.5,0){\vector(1,0){4}}
                 \put(10,4){\vector(0,-1){2}}
                       \put(13,3){\makebox(0,0)[c]{$_{_{\eta^{DT}_{_{_{\overline{X}_i}}}}}$}}
 \put(30,0){\makebox(0,0)[c]{$_{\mb{Q}_{DT}\mb{S}_T(X)_i}$}}
                       \put(38,0){\vector(1,0){2}}
                             \put(39,1){\makebox(0,0)[c]{$_{\lambda^2_i}$}}
                       \put(30,4){\vector(0,-1){2}}
                       \put(31,3){\makebox(0,0)[c]{$_{\phi_i}$}}

\put(53,0){\makebox(0,0)[c]{$_{\mr{Hom}_R(T,A_{X_i})\otimes_S T\oplus P_{A_{X_{i+1}}}\otimes_R DT\otimes_S T}$}}
                            \put(65,0){\vector(1,0){2}}
                 \put(54,4){\line(0,-1){2}}
                 \put(54.2,4){\line(0,-1){2}}
 \put(68,0){\makebox(0,0)[c]{$_0$}}

 \put(0,6){\makebox(0,0)[c]{$_0$}}
                             \put(1,6){\vector(1,0){2}}
 \put(10,6){\makebox(0,0)[c]{$_{\overline{X}_i}$}}
                             \put(17,7){\makebox(0,0)[c]{$_{(-q_i,\overline{u}_{_{X_i}},0)}$}}
                             \put(14,6){\vector(1,0){7}}
 \put(30,6){\makebox(0,0)[c]{$_{X_i\oplus P_{A_{X_i}}\oplus P_{A_{X_{i+1}}}\otimes_R DR}$}}
                             \put(42,8){\makebox(0,0)[c]{\tiny $_{\32jz{u_{_{X_i}}}{0}{p_i}{0}{0}{1}}$}}
                             \put(39,6){\vector(1,0){6}}
 \put(54,6){\makebox(0,0)[c]{$_{A_{X_i}\oplus P_{A_{X_{i+1}}}\otimes DR}$}}
                             \put(64,6){\vector(1,0){2}}
 \put(68,6){\makebox(0,0)[c]{$_0$}}

\end{picture}

\bskip\

Then, since $\overline{X}_i\in\mc{B}$  implies that $\eta_{\overline{X}_i}^{DT}$ is an isomorphism, we obtain that $\phi_i$ is also an isomorphism from the above commutative diagram.


\vskip 10pt
\Exesubsub\label{phi-up}
\hangafter 1\hangindent 40pt
{\it The upper row in the diagram $(\ast)$ is exact
}

\vskip 10pt

In fact, the pullback of $p_i: P_{A_{X_i}}\to A_{X_i}$ and $u_{_{X_i}}: X_i\to A_{X_i}$ in \ref{QSX} gives an exact sequence

\vskip 5pt
\centerline{$0\to \overline{X}_i\stackrel{(-q_i,\overline{u}_{X_i})}{\longrightarrow} X_i\oplus P_{A_{X_i}}\stackrel{{\tiny _{\ljz{u_{_{X_i}}}{p_i}}}}{\longrightarrow} A_{X_i}\to 0$,}
\vskip 5pt

\noindent since $p_i$ is surjective. The direct sum of the above exact sequence and the trivial exact sequence $0\to 0\to P_{A_{X_{i+1}}}\stackrel{1}{\longrightarrow}P_{A_{X_{i+1}}}\to 0$ gives us the exact sequence in the upper row in the diagram $(\ast)$.


\vskip 10pt
\Exesubsub\label{phi-bottom}
\hangafter 1\hangindent 40pt
{\it The bottom row in the diagram $(\ast)$
}

\vskip 10pt

Note that we have the following exact sequence in \ref{QSX}

\vskip 5pt
\centerline{$0\to \mb{Q}_{DT}\mb{S}_T(X)\stackrel{\lambda}{\longrightarrow} \mr{R}_{DT}(\mf{L}^P_{A_X})\stackrel{r(s^P)}{\longrightarrow} \mr{\hat{H}om}(DT, \mb{S}_T(X))\to 0$,}
\vskip 5pt

\noindent and that

\vskip 5pt
\hangafter 1\hangindent 20pt
$\mr{R}_{DT}(\mf{L}^P_{A_X})_i=\mr{Hom}_S(DT, (\mf{L}^P_{A_X})_i)\oplus  (\mf{L}^P_{A_X})_{i+1}\otimes_S T$\\
 $= \mr{Hom}_S(DT, \mr{Hom}_R(T,A_{X_{i-1}})\oplus P_{A_{X_i}}\otimes_R DT)\oplus  (\mr{Hom}_R(T,A_{X_i})\oplus  P_{A_{X_{i+1}}}\otimes_R DT)\otimes_S T$.
\vskip 5pt

\noindent So we have the following pullback diagram, for some homomorphisms $a_i$,

\bskip\

 \setlength{\unitlength}{0.09in}
 \begin{picture}(50,7)

 \put(0,0){\makebox(0,0)[c]{$_0$}}
                             \put(1,0){\vector(1,0){2}}

 \put(9,0){\makebox(0,0)[c]{$_{\mr{Hom}_R(DT, \overline{X}_i\otimes_R DT)}$}}
                             \put(17,1){\makebox(0,0)[c]{$_{b_i}$}}
                             \put(16,0){\vector(1,0){2}}
                             \put(9,4){\line(0,-1){2}}
                             \put(9.2,4){\line(0,-1){2}}
 \put(33,0){\makebox(0,0)[c]{$_{_{\mr{Hom}_S(DT, \mr{Hom}_R(T,A_{X_{i-1}})\oplus P_{A_{X_i}}\otimes_R DT)}}$}}
                       \put(45,0){\vector(1,0){2}}
                             \put(46,1){\makebox(0,0)[c]{$_{r^1_i}$}}
                       \put(30,4){\vector(0,-1){2}}
                       \put(31,3){\makebox(0,0)[c]{$_{\lambda^1_i}$}}

\put(55,0){\makebox(0,0)[c]{$_{\mr{Hom}_S(DT,\mb{S}_T(X)_i)}$}}
                            \put(66,0){\vector(1,0){2}}
                 \put(55,4){\vector(0,-1){2}}
                 \put(57,3){\makebox(0,0)[c]{$_{-r^2_i}$}}
 \put(70,0){\makebox(0,0)[c]{$_0$}}

 \put(0,6){\makebox(0,0)[c]{$_0$}}
                             \put(1,6){\vector(1,0){2}}
 \put(9,6){\makebox(0,0)[c]{$_{\mr{Hom}_S(DT, \overline{X}_i\otimes_R DT)}$}}
                             \put(17,7){\makebox(0,0)[c]{$_{a_i}$}}
                             \put(16,6){\vector(1,0){2}}
 \put(30,6){\makebox(0,0)[c]{$_{\mb{Q}_{DT}\mb{S}_T(X)_i}$}}
                             \put(43,7){\makebox(0,0)[c]{$_{\lambda^2_i}$}}
                             \put(42,6){\vector(1,0){2}}
 \put(55,6){\makebox(0,0)[c]{$_{(\mr{Hom}_R(T,A_{X_{i-1}})\oplus P_{A_{X_i}}\otimes_R DT)\otimes_S T}$}}
                             \put(66,6){\vector(1,0){2}}
 \put(70,6){\makebox(0,0)[c]{$_0$}}

\end{picture}
\bskip\

\noindent where $r^1_i, r^2_i$ and $\lambda^1_i, \lambda^2_i$ are the components of the homomorphisms $r(s^P)_i$ and $\lambda_i$ respectively, and $b_i=\mr{Hom}_S(DT,t_i)$ with

\vskip 5pt
\centerline{$t_i=$ {\footnotesize $(-(q_i\otimes DT)\circ \theta_{l_i},\overline{u}_{X_i}\otimes_R DT): \overline{X}_i\otimes_RDT\to \mr{Hom}_R(T,A_{X_{i-1}})\oplus P_{A_{X_i}}\otimes_R DT$}}
\vskip 5pt

\noindent is given in \ref{QSX}.

Note that $r^1_i$ is surjective as we indicate in \ref{ker} in the general case, so the upper row is exact. Thus we get the bottom exact sequence in the diagram $(\ast)$.


\vskip 10pt
\Exesubsub\label{phi-comm}
\hangafter 1\hangindent 40pt
{\it The  diagram $(\ast)$ is commutative
}

\vskip 10pt

At first, it is easy to see that
 the right part of the diagram $(\ast)$ is commutative from the construction of the morphism $\phi$ in 3.4.1.

As to  the left part of the diagram $(\ast)$, we first show that the following equality of  compositions

\hskip 120pt{$\eta^{DT}_{_{\overline{X}_i}}\circ a_i \circ \lambda^1_i=(-q_i,\overline{u}_{X_i},0)\circ \phi_i\circ \lambda^1_i$. \hfill $(\dag_1)$}

Indeed, we have that

\vskip 5pt
\hskip 20pt
\hangafter 1\hangindent 50pt
$\eta^{DT}_{_{\overline{X}_i}}\circ a_i \circ \lambda^1_i$\\
$=\eta^{DT}_{_{\overline{X}_i}}\circ b_i $ \hfill (by the commutative diagram in \ref{phi-bottom})\\
$= \eta^{DT}_{_{\overline{X}_i}}\circ \mr{Hom}_S(DT,t_i)$ \hfill (since $b_i=\mr{Hom}_S(DT,t_i)$)\\
$=\mb{\Gamma}^{DT}(1_{\overline{X}_i\otimes_RDT}) \circ  \mr{Hom}_S(DT,t)$\\
$=\mb{\Gamma}^{DT}(1_{\overline{X}_i\otimes_RDT}\circ t)$ \hfill (by the naturality  of $\mb{\Gamma}^{DT}$)\\
$=\mb{\Gamma}^{DT}(t)$\\
$=\mb{\Gamma}^{DT}(${\footnotesize $(-(q_i\otimes_R DT)\circ \theta_{l_i},\overline{u}_{X_i}\otimes_R DT)$}$)$ \hfill (since  {\footnotesize $t_i=(-(q_i\otimes DT)\circ \theta_{l_i},\overline{u}_{X_i}\otimes_R DT)$})

\vskip 5pt

\noindent and we also have that

\vskip 5pt
\hskip 20pt
\hangafter 1\hangindent 50pt
$(-q_i,\overline{u}_{X_i},0)\circ \phi_i\circ \lambda^1_i$\\
$=(-q_i,\overline{u}_{X_i},0)\circ \xi_i^a$ \hfill (since $\xi_i=(\xi_i^a,\xi_i^b)=\phi\circ \lambda$)\\
$=(-q_i,\overline{u}_{X_i},0)\circ ${\footnotesize $\32jz{\xi^a_{11}}{0}{0}{\xi^a_{22}}{0}{0}$} \hfill (by \ref{xi-a})\\
$=(-q_i\circ \xi^a_{11}, \overline{u}_{X_i}\circ \xi^a_{22}) $.\\
\vskip 5pt

\noindent But $\xi^a_{11}=\mb{\Gamma}^{DT}(\theta_{l_i})$ and $\xi^a_{22}=\mb{\Gamma}^{DT}(1_{(P_{A_{X_i}}\otimes_R DT)})$ by the construction in \ref{xi-a}, we obtain that

\vskip 5pt
\centerline{$q_i\circ \xi^a_{11}=q_i\circ \mb{\Gamma}^{DT}(\theta_{l_i})=\mb{\Gamma}^{DT}(q_i\otimes_RDT\circ \theta_{l_i})$}
\vskip 5pt

\noindent and that

\vskip 5pt
\centerline{$\overline{u}_{X_i}\circ \xi^a_{22}=\overline{u}_{X_i}\circ \mb{\Gamma}^{DT}(1_{(P_{A_{X_i}}\otimes_R DT)})=\mb{\Gamma}^{DT}(\overline{u}_{X_i}\otimes_RDT\circ 1_{(P_{A_{X_i}}\otimes_R DT)})=\mb{\Gamma}^{DT}(\overline{u}_{X_i}\otimes_RDT)$}
\vskip 5pt

Hence, we see that the equality $(\dag_1)$ holds.

\vskip 5pt
Since that $a_i\circ\lambda^2_i=0$ and that

\centerline{$(-q_i,\overline{u}_{X_i},0)\circ \phi_i\circ \lambda^2_i=(-q_i,\overline{u}_{X_i},0)\circ \xi^b_i=0$,}

\noindent we also get that

\hskip 120pt{$\eta^{DT}_{_{\overline{X}_i}}\circ a_i \circ \lambda^2_i=(-q_i,\overline{u}_{X_i},0)\circ \phi_i\circ \lambda^2_i$. \hfill $(\dag_2)$}

\vskip 5pt
\noindent Now, from the property of  the pullback in \ref{phi-bottom}, we know that the two equalities $(\dag_1)$ and  $(\dag_2)$ together imply that

\centerline{$\eta^{DT}_{_{\overline{X}_i}}\circ a_i=(-q_i,\overline{u}_{X_i},0)\circ \phi_i$.}

\noindent Thus, the left part of the diagram is also commutative.

\vskip 5pt


\vskip 10pt
\Exesub\label{phi-nat}
\hangafter 1\hangindent 40pt
{\it The isomorphism $\phi: X\oplus \mr{L}_R(P^+_{A_X})\to \mb{Q}_{DT}\mb{S}_T(X)$ is natural on $X$.
}

\vskip 10pt

For any $X,Y\in\mc{RC}(R)$ and $f\in\mr{Hom}_{\mc{RC}(R)}(X,Y)$, it is regular to show that the following diagram is commutative, for some natural homomorphisms $f\oplus \mr{L}_R(P^+_{A_f})$ and $\mb{Q}_{DT}\mb{S}_T(f)$.

\vskip 15pt

 \setlength{\unitlength}{0.09in}
 \begin{picture}(50,7)

  \put(20,0){\makebox(0,0)[c]{$Y\oplus \mr{L}_R(P^+_{A_Y})$}}
                             \put(29,1.5){\makebox(0,0)[c]{$_{\phi_{_Y}}$}}
                             \put(27,0){\vector(1,0){4}}
 \put(39,0){\makebox(0,0)[c]{$\mb{Q}_{DT}\mb{S}_T(X)$.}}

                 \put(20,4){\vector(0,-1){2}}
                       \put(16,3){\makebox(0,0)[c]{$_{f\oplus \mr{L}_R(P^+_{A_f})}$}}
                 \put(38,4){\vector(0,-1){2}}
                       \put(42,3){\makebox(0,0)[c]{$_{\mb{Q}_{DT}\mb{S}_T(f)}$}}

 \put(20,6){\makebox(0,0)[c]{$X\oplus \mr{L}_R(P^+_{A_X})$}}
                             \put(29,7.5){\makebox(0,0)[c]{$_{\phi_{_X}}$}}
                             \put(27,6){\vector(1,0){4}}
 \put(39,6){\makebox(0,0)[c]{$\mb{Q}_{DT}\mb{S}_T(X)$}}

\end{picture}
\vskip 10pt

\noindent Thus, the isomorphism $\phi$ is natural on $X$. This means that $\mb{Q}_{DT}\mb{S}_T\simeq 1_{\underline{\mc{RC}}(R)}$ naturally.


\zskip\

\Exe\label{SQ=1}
{\bf The isomorphism  $\mb{S}_T\mb{Q}_{DT}\simeq 1_{\underline{\mc{RC}}(S)}$}
%
%
\bskip\

Dually to the proof of \ref{QS=1}, one can show that  $\mb{S}_T\mb{Q}_{DT}\simeq 1_{\underline{\mc{RC}}(S)}$ naturally.

Namely, for an object $Y\in\mc{RC}(S)$, one uses that $(\mc{G},\mc{K})$ is a complete hereditary cotorsion pair in $\mr{mod}S$ to obtain exact sequences $0\to K_{Y_i}\to G_{Y_i}\to Y_i\to 0$, for each $i$. Then take an injective $S$-module $I_{G_{Y_i}}$ and a monomorphism $j: G_{Y_i}\to I_{G_{Y_i}}$, one can show that there is a natural isomorphism $\mb{S}_T\mb{Q}_{DT}(Y)\to Y\oplus \mr{R}_{DS}(I^-_Y)$, where $I^-_Y=\{I_{G_{Y_{i-1}}}\}$ and $\mr{R}_{DS}(I^-_Y)$ is a projective object in $\mc{RC}(R)$ by Remark (2) in \ref{rFct}. And then one gets that  $\mb{S}_T\mb{Q}_{DT}\simeq 1_{\underline{\mc{RC}}(S)}$ naturally.


\zskip\

\Exe\label{pfth}
{\bf The last proof of main theorem}
%
%
\zskip\
Recall that $[1]$ is an automorphism of repetitive categories, where $(X[1])_i=X_{i-1}$ for an object in a repetitive category.

Define $\mb{F}_T:=[-1]\mb{S}_T: \underline{\mc{RC}}(R)\to \underline{\mc{RC}}(S)$ and $\mb{G}_T:=\mb{Q}_{DT}[1]: \underline{\mc{RC}}(S)\to \underline{\mc{RC}}(R)$. Then we have that $\mb{F}_T\mb{G}_T\simeq 1_{\underline{\mc{RC}}(S)}$ naturally and that  $\mb{G}_T\mb{F}_T\simeq 1_{\underline{\mc{RC}}(R)}$ naturally. So that $\mb{F}_T$ and $\mb{G}_T$ gives a repetitive equivalence between $R$ and $S$.

It is easy to check that $\mb{F}_T|_{\mc{A}}\simeq \mr{Hom}_R(T,-)$ and that  $\mb{G}_T|_{\mc{G}}\simeq -\otimes_ST$ from the definitions of two functors. Now the proof of the theorem is completed.


\zskip\

\Exe\label{pfth}
{\bf The  proof of the proposition in the introduction}
%
%
\bskip\
 Assume that the equivalence is given by the functor $F:\underline{\mc{RC}}(R)\to\underline{\mc{RC}}(S)$. By assumptions, $F$ restricts to an equivalence $\mc{A}\to \mc{G}$. Note that $\mc{G}$ is resolving and $S\in\mc{G}$. Let $T=F^{-1}(S)$. Then $T\in\mc{A}$. By the triangle equivalence, we have that, for any $A\in\mc{A}$, \\
\centerline{$\mr{Ext}_R^i(T,A)\simeq\mr{Hom}_{\underline{\mc{RC}}(R)}(T,\Sigma^iA)\simeq \mr{Hom}_{\underline{\mc{RC}}(S)}(S,\Sigma^iF(A))\simeq \mr{Ext}_S^i(S,F(A))$,}
where $\Sigma$ is the translator funtor in repetitive categories. In particular, we obtain that $\mr{Hom}_R(T,A)\simeq \mr{Hom}_S(S,F(A))\simeq F(A)$ and that $\mr{Ext}_R^i(T,A)=0$ for all $i>0$.
It follows that  $S\simeq\mr{End}(T_R)$ and $\mr{Ext}_R^i(T,T)=0$ for all $i>0$.
Note that $\mc{A}$ is coresolving and $DR\in\mc{A}$, so we also have that $F(DR)\simeq\mr{Hom}_R(T,DR)\simeq DT$.   Thus,
we get that  \\
{\footnotesize\centerline{$\mr{Ext}_S^i({_ST},{_ST})\simeq\mr{Ext}_S^i(DT,DT)\simeq\mr{Hom}_{\underline{\mc{RC}}(S)}(DT,\Sigma^iDT)$\hskip 110pt} \\
\centerline{$\simeq \mr{Hom}_{\underline{\mc{RC}}(S)}(F(DR),\Sigma^iF(DR))\simeq\mr{Hom}_{\underline{\mc{RC}}(R)}(DR,\Sigma^iDR)\simeq \mr{Ext}_R^i(DR,DR)$.}} It follows that $\mr{End}(_ST)\simeq R$ and that $\mr{Ext}_S^i({_ST},{_ST})=0$ for all $i>0$. Thus, $T$ is a Wakamatsu-tilting module.

By assumption, $F|_{\mc{A}}\simeq\mr{Hom}_R(T,-)$ gives the equivalence $\mc{A}\to \mc{G}$. It follows that $F^{-1}|_{\mc{G}}\simeq -\otimes_ST$ by the unique of the adjoint. Note that $\mr{Hom}_R(T,-)$ and $-\otimes_ST$ are exact functors respectively in $\mc{A}$ and ${\mc{G}}$, since $F$ is a triangle functor. As $\mc{A}$ is coresolving, for any $A\in\mc{A}$, the exact sequence $0\to A\to I\to A'\to 0$ with $I\in\mr{inj}R$ is a sequence in $\mc{A}$. Applying the exact functor $\mr{Hom}_R(T,-)$, we obtain that $\mr{Ext}_R^1(T,A)=0$. It follows that $T\in\mc{A}\bigcap\mr{KerExt}^1_R(-,\mc{A})$, i.e, $T$ is Ext-projective in $\mc{A}$. Dually, we have also that $\mr{Tor}^S_1(X,T)=0$ for any $X\in\mc{G}$. In particular, for any $A\in\mc{A}$, suppose that $A=X\otimes_ST$ for some $X\in\mc{G}$ and take an exact sequence $0\to X'\to P\to X\to 0$ with $P\in\mr{proj}R$, then the sequence is in $\mc{G}$ since $\mc{G}$ is resolving, and hence there is an induced exact sequence $0\to  X'\otimes_ST\to P\otimes_ST\to X\otimes_ST\to 0$, since $-\otimes_ST$ is exact in $\mc{G}$. The last sequence gives an exact sequence $0\to A''\to T_A\to A\to 0$ with $T_A= P\otimes_ST\in\mr{add}_RT$ and $A''=X'\to P\in\mc{A}$. It follows that $T$ is an Ext-projective generator in $\mc{A}$. Now applying Corollary \ref{c-gwt}, we conclude that $T$ is a good Wakamatsu-tilting module.
%
%


%
\vskip 30pt
\section{\large Examples}


\Exe\label{t-ct}
{\bf Tilting modules and cotilting modules }

\vskip 10pt
Let $R$ ba an artin algebra. Recall that an $R$-module $T$ is tilting provided the following three conditions are satisfied:

\vskip 5pt
\hangafter 1\hangindent 40pt
(1) The projective dimension of $T$ is finite;

(2) $\mr{Ext}_R^i(T,T)=0$ for all $i>0$;

\hangafter 1\hangindent 40pt
(3) There is an exact sequence $0\to R\to T_0\to\cdots\to T_n\to 0$ for some integer $n$, where each $T_i\in\mr{add}_RT$.

\vskip 5pt
\noindent
Dually, an $R$-module $T$ is cotilting provided the following three conditions are satisfied:

\vskip 5pt
\hangafter 1\hangindent 40pt
(1) The injective dimension of $T$ is finite;

(2) $\mr{Ext}_R^i(T,T)=0$ for all $i>0$;

\hangafter 1\hangindent 40pt
(3) There is an exact sequence $0\to T_n\to\cdots\to T_0\to DR\to 0$ for some integer $n$, where each $T_i\in\mr{add}_RT$.

\vskip 5pt
It is known that an $R$-module $T$ is a tilting module if and only if $DT$ is a cotilting left $R$-module if and only if $DT$ is a cotilting $S$-module. Note also that both tilting modules and cotilting modules are Wakamatsu-tilting modules.

We need the following well-known results on tilting modules and cotilting modules.

\vskip 10pt
\noindent\hangafter 1\hangindent 38pt
{\bf Proposition} {\it %
$(1)$ If $T$ is a tilting module, then the cortorsion pair $(\mr{KerExt}_R^1(-,{_T\mc{X}}),{_T\mc{X}})$ is complete.

\hangafter 1\hangindent 40pt
$(2)$ If $T$ is a cotilting module, then the cortorsion pair $(\mc{X}_T,\mr{KerExt}_R^1(\mc{X}_T,-))$ is complete.
}

\vskip 5pt
\noindent
\Pf. (2)  follows from \cite[Section 5]{AR}, and (1) is just the dual of (2).
\hfill $\Box$

\zskip\

\Exesub\label{tm}
{\it Tilting modules are  good Wakamatsu-tilting}

\vskip 10pt

Assume $T_R$ is a tilting module of finite projective dimension. Let $S=\mr{End}(T_R)$. Then $_ST_R$ is a good Wakamatsu-tilting module. Hence there is an equivalence between repetitive categories $\underline{\mc{RC}}(R)$ and $\underline{\mc{RC}}(S)$.

Indeed, if $_ST_R$ is a tilting module of finite projective dimension, then  $T$ is Wakamatsu-tilting and $_RDT_S$ is a cotilting module of finite injective dimension. By Proposition \ref{cotp}, we obtain that  bimodules $_ST_R$ and $_RDT_S$ represent a cotorsion pair counter equivalence between the compete hereditary cotorsion pair $(\mr{KerExt}_R^1(-,{_T\mc{X}}),{_T\mc{X}})$ in $\mr{mod}R$ and the compete hereditary cotorsion pair  $(\mc{X}_{DT},\mr{KerExt}_S^1(\mc{X}_{DT},-))$ in $\mr{mod}S$. It follows from the definition that $_ST_R$ is a good Wakamatsu-tilting bimodule.

\zskip\

\Exesub\label{tm}
{\it Cotilting modules are  good Wakamatsu-tilting}

\vskip 10pt

Assume now $T_R$ is a cotilting module of finite injective dimension with $S=\mr{End}(T_R)$. Then $_ST_R$ is also a good Wakamatsu-tilting module. Hence there is an equivalence between repetitive categories $\underline{\mc{RC}}(R)$ and $\underline{\mc{RC}}(S)$.

Indeed, dually to \ref{tm}, if $_ST_R$ is a cotilting module of finite injective dimension, then $_RDT_S$ is a tilting module of finite projective dimension. By Proposition \ref{cotp} again, we obtain that  bimodules $_ST_R$ and $_RDT_S$ represent a cotorsion pair counter equivalence between the compete hereditary cotorsion pair $(\mc{X}_{T},\mr{KerExt}_S^1(\mc{X}_{T},-))$ in $\mr{mod}R$ and the compete hereditary cotorsion pair $(\mr{KerExt}_R^1(-,{_{DT}\mc{X}}),{_{DT}\mc{X}})$ in $\mr{mod}S$.  It follows from the definition that $_ST_R$ is also a good Wakamatsu-tilting bimodule.

\zskip\

\Exe\label{fpt}
{\bf Wakamatsu-tilting modules of finite  type}

\vskip 10pt

\Exesub\label{wfpt}
%
We say  that a  Wakamatsu-tilting $R$-module $T$ is of finite  type provided that either the subcategory $\mr{KerExt}_R^1(-,{_T\mc{X}})$ or  the subcategory $\mr{KerExt}_R^1(\mc{X}_T,-)$ is of finite representation type.
In particular, if $R$ is an algebra of finite representation type, then each subcategory of $\mr{mod}R$ is of finite representation type, and hence every Wakamatsu-tilting module in $\mr{mod}R$ is of finite  type.

We note that, if $T$ is a  Wakamatsu-tilting $R$-module of finite  type with $S=\mr{End}(T_R)$, then $DT$ is a Wakamatsu-tilting $S$-module of finite type. This is just followed from the equivalences in Proposition \ref{cotp}.

\vskip 10pt
\noindent
{\bf Proposition}\ \ %
{\it A  Wakamatsu-tilting module of finite  type is always a good  Wakamatsu-tilting module. In particular,  every Wakamatsu-tilting module over an  algebra of finite representation type is good.
}

\vskip 5pt
\noindent
\Pf. Let $T$ be a  Wakamatsu-tilting $R$-module of finite  type with $S=\mr{End}(T_R)$. Assume first that  the subcategory $\mr{KerExt}_R^1(\mc{X}_T,-)$ is of finite representation type. Then the hereditary cotorsion pair $(\mc{X}_{T},\mr{KerExt}_S^1(\mc{X}_{T},-))$ in $\mr{mod}R$ is complete. Moreover, by the equivalence in Proposition \ref{cotp} (3), the subcategory $\mr{KerExt}_R^1(-,{_{DT}\mc{X}})$ is also  of finite representation type. Thus, the hereditary cotorsion pair $(\mr{KerExt}_R^1(-,{_{DT}\mc{X}}),{_{DT}\mc{X}})$ in $\mr{mod}S$ is also complete. It follows that  bimodules $_ST_R$ and $_RDT_S$ represent a cotorsion pair counter equivalence between the compete hereditary cotorsion pair $(\mc{X}_{T},\mr{KerExt}_S^1(\mc{X}_{T},-))$ in $\mr{mod}R$ and the compete hereditary cotorsion pair $(\mr{KerExt}_R^1(-,{_{DT}\mc{X}}),{_{DT}\mc{X}})$ in $\mr{mod}S$.  Similarly, in case that $\mr{KerExt}_R^1(-,{_T\mc{X}})$ is of finite representation type, we have that  bimodules $_ST_R$ and $_RDT_S$ represent a cotorsion pair counter equivalence between the compete hereditary cotorsion pair $(\mr{KerExt}_R^1(-,{_T\mc{X}}),{_T\mc{X}})$ in $\mr{mod}R$ and the compete hereditary cotorsion pair  $(\mc{X}_{DT},\mr{KerExt}_S^1(\mc{X}_{DT},-))$ in $\mr{mod}S$. Altogether, we see that $T$ is a good Wakamatsu-tilting module in either case.
\hfill $\Box$

\vskip 10pt

\Exesub\label{wfpt}
%
Two trivial examples of Wakamatsu-tilting modules of finite  type over an algebra $R$ is the module $R$ and the module $DR$. In the first case, the subcategory $\mr{KerExt}_R^1(-,{_T\mc{X}})=\mr{proj}R$ is of finite representation type, while the subcategory $\mr{KerExt}_R^1(\mc{X}_T,-)=\mr{inj}R$ is of finite representation type in the second case.

The following is an example of Wakamatsu-tilting modules of finite  type over an algebra of infinite representation type.

\vskip 10pt
\noindent{\bf Example} Let $R$ be the bounded quiver algebra given by the following quiver over a field with the relation given by $\mr{rad}^2R=0$.

    $$\xymatrix{
        1 \ar@/^/[r]^{\alpha} & \ar@/^/[l]^{\beta} 2  \ar[r] \ar@{.>}[r]^{\gamma} & 3 \ar[r] \ar@{.>}[r]^{\delta} & 4
         \ar@/^/[r]^{\epsilon} & \ar@/^/[l]^{\varepsilon} 5  \ar[r] \ar@{.>}[r]^{\zeta} & 6 \ar@<0.5ex>[r]^{\eta} \ar[r]_{\theta} & 7}$$

The following is the AR-quiver of the algebra.
%
{\scriptsize
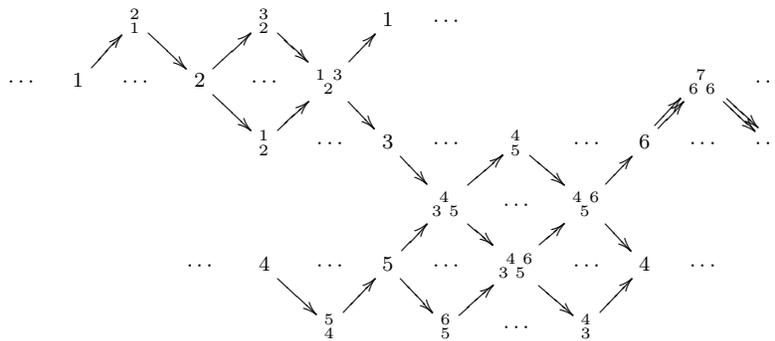
\noindent    $$\xymatrix@C=1em@R=8pt{
        &    & {^2_1}\ar[dr] &     & {^3_2}\ar[dr]&       & 1  & \cdots     &   &    &               &     &     \\
   \cdots & 1\ar[ur] & \cdots      &2\ar[dr]\ar[ur]&\cdots  & {^{1\ 3}_{\ 2\ }}\ar[dr]\ar[ur]& &   &     &     &     &  {^{\ 7\ }_{6\ 6}}  \ar@<0.5ex>[dr]\ar[dr] & \cdots\\
       &      &             &        &{^1_2}\ar[ur]&\cdots&3\ar[dr]&\cdots&{^4_5}\ar[dr]&\cdots&6\ar@<0.5ex>[ur]\ar[ur] &  \cdots  &\cdots \\
       &       &             &       &             &       &       &{^{\ 4\ }_{3\ 5}}\ar[ur]\ar[dr]&\cdots &{^{4\ 6}_{\ 5\ }}\ar[ur]\ar[dr]& & &  \\
       &      &       &\cdots   & 4\ar[dr] &\cdots &5\ar[ur]\ar[dr]&\cdots &{^{\ 4\ 6}_{3\ 5\ }}\ar[ur]\ar[dr]&\cdots& 4 &\cdots &  \\
       &      &       &         &   &{^5_4}\ar[ur] &  & {^6_5}\ar[ur] & \cdots & {^4_3} \ar[ur] & & & }$$
}

The algebra is of infinite representation type. Over this algebra, we have a  Wakamatsu-tilting module of finite  type (and hence, a good Wakamatsu-tilting module)

\vskip 5pt
\centerline{$T={\footnotesize \begin{array}{c}   2\\ 1\end{array}} \oplus {\footnotesize \begin{array}{c} 1\ 3\\2 \end{array}} \oplus {\footnotesize \ 3\ } \oplus {\footnotesize \begin{array}{c} 4\\ 3\ 5 \end{array}}\oplus {\footnotesize \begin{array}{c} 5\\4 \end{array}} \oplus {\footnotesize \begin{array}{c} 6\\5 \end{array}} \oplus {\footnotesize \begin{array}{c} 7\\ 6\ 6 \end{array}} $.}
%

\vskip 5pt
Indeed, one can check that the subcategory $\mr{KerExt}_R^1(-,{_T\mc{X}})$ is of finite representation type, while the subcategory ${_T\mc{X}}$ is of infinite representation type.

%

%
%


\vskip 30pt

{\small

}

%
%

\end{document}